\documentclass[smallextended,referee,envcountsect]{svjour3}
\smartqed
\usepackage{graphicx}
\usepackage[hyperfootnotes=false]{hyperref}
\hypersetup{colorlinks=true,citecolor=blue,linkcolor=blue,urlcolor=blue}
\usepackage{amsmath,amsxtra,amssymb,latexsym,amsfonts,amscd}
\usepackage{multicol,color}
\usepackage{float}
\usepackage{soul}
\usepackage{graphicx}
\usepackage{caption}
\usepackage{subcaption}
\usepackage{pstricks-add}
\usepackage{pgf,tikz}
\usepackage{mathrsfs}
\usetikzlibrary{arrows}
\usepackage{amssymb}
\usepackage{theorem}
\usepackage{fancyhdr}
\usepackage{tikz}
\usepackage{enumerate}
\usepackage[margin=1.1in]{geometry}
\def\disp{\displaystyle}
\def\e{\epsilon}

\def\DD{\Delta}
\def\lm{\lambda}

\def\({\left(}
\def\){\right)}
\def\[{\left[}
\def\]{\right]}
\def\n{\left \|}
\def\en{\right \|}

\def\ox{\bar{x}}

\def\ou{\bar{u}}

\def\gg{\gamma}
\def\vt{\vartheta}
\def\dn{\downarrow}

\def\la{\langle}
\def\ra{\rangle}

\def\e{\varepsilon}

\def\h{\hfill\Box}
\def\R{\mathbb{R}}
\def\N{\mathbb{N}}

\def\dn{\downarrow}

\def\vph{\varphi}
\def\emp{\emptyset}

\def\lm{\lambda}

\def\gg{\gamma}

\def\DD{\Delta}
\def\al{\alpha}

\def\be{\beta}

\def\ph{\varphi}

\def\N{I\!\!N}
\def\th{\theta}

\newcounter{lk}

\theoremstyle{plain}{\theorembodyfont{\rmfamily}
}
\theoremstyle{plain}{\theorembodyfont{\rmfamily}
}
\theoremstyle{plain}{\theorembodyfont{\rmfamily}
}
\theoremstyle{plain}{\theorembodyfont{\rmfamily}

\theoremstyle{plain}{\theorembodyfont{\rmfamily}

\def\eq{\begin{equation}}
\def\eeq{\end{equation}}
\journalname{JOTA}
\usepackage{atbegshi}

\begin{document}

\title{Optimal Control of Sweeping Processes in Robotics and Traffic Flow Models}

\author{Giovanni Colombo$\;{\bf\cdot}\;$Boris Mordukhovich$\;{\bf\cdot}\;$Dao Nguyen}

\institute{Giovanni Colombo \at
Dipartimento di Matematica Pura e Applicata, Universit$\grave{\mbox{a}}$ di Padova\\
via Trieste 63, 35121 Padova, Italy\\
colombo@math.unipd.it
\and
Boris S. Mordukhovich,  Corresponding author \at
Department of Mathematics, Wayne State University\\
Detroit, Michigan, USA\\
boris@math.wayne.edu
\and
Dao Nguyen \at
Department of Mathematics, Wayne State University\\
Detroit, Michigan, USA\\
dao.nguyen2@wayne.edu}

\date{Received: date / Accepted: date}

\maketitle\vspace*{-0.2in}

\begin{abstract}
The paper is mostly devoted to applications of a novel optimal control theory for perturbed sweeping/Moreau processes to two practical dynamical models. The first model addresses mobile robot dynamics with obstacles, and the second one concerns control and optimization of traffic flows. Describing these models as controlled sweeping processes with pointwise/hard control and state constraints and applying new necessary optimality conditions for such systems allow us to develop efficient procedures to solve naturally formulated optimal control problems for the models under consideration and completely calculate optimal solutions in particular situations.
\end{abstract}\vspace*{-0.1in}
\vspace*{-0.1in}

\vspace*{-0.25in}

\section{Introduction}\label{Introduction}\vspace*{-0.15in}

Sweeping process models were introduced by Jean-Jacques Moreau in the 1970s to describe dynamical processes arising in elastoplasticity and related mechanical areas; see \cite{mor_frict}. Such models were given in the form of discontinuous differential inclusions governed by the normal cone mappings to nicely moving convex sets. It has been well realized in the sweeping process theory that the Cauchy problem for the basic Moreau's sweeping process and its slightly nonconvex extensions admits unique solutions; see, e.g., \cite{CT}. This therefore excludes any possible optimization of sweeping differential inclusions and strikingly distinguishes them from the well-developed optimal control theory for their Lipschitzian counterparts. On the other hand, existence and uniqueness results for sweeping trajectories provide a convenient framework for handling simulation and related issues in various applications to mechanics, hysteresis, economics, robotics, electronics, etc.; see, e.g., \cite{adly,brog,HB,jv,mv2} among more recent publications with the references therein.\vspace*{-0.05in}

To the best of our knowledge, first control problems associated with sweeping processes and first topics to investigate were related to the existence and relaxation of optimal solutions to sweeping differential inclusions with controls in additive perturbations as developed and discussed in \cite{et}. Starting with \cite{chhm}, serious attention has been drawn to optimal control problems for sweeping processes with control actions entering moving sets and deriving necessary optimality conditions in various state-constrained optimal control problems that appear in this way for discontinuous sweeping differential inclusions; see \cite{chhm1,cm2,cm3,cm4,hm18}. Advanced necessary optimality condition for control systems governed by sweeping processes with constrained controls in additive perturbations have been recently derived in \cite{bk,ao,ac,pfs,cmn}.\vspace*{-0.05in}

In this paper we present new applications of the most recent necessary optimality conditions obtained in our paper \cite{cmn} to two classes of practical models. The first one is taken from the area of robotics, while the second model concerns pedestrian traffic flows. Dynamics in these models can be formalized as a perturbed sweeping process. Inserting constrained control actions into a perturbation force and selecting a practically motivated cost functional allow us to describe the corresponding controlled dynamical systems in the form of optimal control problems studied in \cite{cmn}. Then we apply the necessary optimality condition from \cite{cmn} to the obtained control problems and express them entirely in terms of the given data. This brings us to precise relationships for computing optimal solutions in some major situations, which are discussed in detail and are illustrated by nontrivial examples.\vspace*{-0.05in}

The rest of the paper is organized as follows. In Section~2 we recall for the reader's convenience the results of \cite{cmn} needed for our subsequent applications. Section~3 is devoted to formulating and solving an optimal control version of the mobile robot model with obstacles that is well recognized in robotics. Section~4 deals with a deterministic continuous-time optimal control version of the pedestrian traffic flow model that belongs to the area of socioeconomics. The concluding Section~5 presents a summary of the major results and discusses some unsolved problems of the future research.\vspace*{-0.05in}

Throughout the paper we use standard notations from variational analysis, control theory, and the applied areas of modeling, which are specified in the corresponding places below. Recall here that, given a matrix $A$, the symbol $A^*$ indicates its transposition/adjoint operator.\vspace*{-0.28in}

\section{Discretization and Necessary Optimality Conditions for Controlled Sweeping Processes}\label{Assumptions}\vspace*{-0.15in}

In this section we formulate the general optimal control problem for a perturbed sweeping process studied in \cite{cmn} and present some major results of that paper needed in the sequel.\vspace*{-0.05in}

Denote by $(P)$ the following optimal control problem:
\begin{eqnarray}\label{cost1}
\mbox{minimize}\quad J[x,u]:=\vph\big(x(T)\big)
\end{eqnarray}
over pairs $(x(\cdot),u(\cdot))$ of measurable controls $u(t)$ and absolutely continuous trajectories $x(t)$ on the fixed time interval $[0,T]$ satisfying the {\em controlled sweeping differential inclusion}
\begin{eqnarray}\label{1.2}
\dot{x}(t)\in-N\big(x(t);C\big)+g\big(x(t),u(t)\big)\;\mbox{ a.e. }\;t\in [0,T],\;x(0):=x_0\in C\subset\R^n,
\end{eqnarray}
subject to the {\em pointwise constraints on control} actions
\begin{eqnarray}\label{cont}
u(t)\in U\subset\R^d\;\mbox{ a.e. }\;t\in[0,T].
\end{eqnarray}
The set $C$ in \eqref{1.2} is a {\em convex polyhedron} given by
\begin{eqnarray}\label{C}
C:=\bigcap_{j=1}^s C^j\;\mbox{ with }\;C^j:=\big\{x\in\R^n\big|\;\la x^j_*,x\ra\le c_j\big\},
\end{eqnarray}
and the normal cone to it in \eqref{1.2} is understood in the classical sense of convex analysis
\begin{eqnarray}\label{1.4}
N\big(x;C\big):=\big\{v\in\R^n\big|\;\la v,y-x\ra\le 0,\;y\in C\big\}\mbox{ if }x\in C\;\mbox{ and }\;N\big(x;C\big):=\emp\;\mbox{ if }\;x\notin C.
\end{eqnarray}
It follows directly from \eqref{1.2} due to the second part of the normal cone definition \eqref{1.4} that we implicitly have the {\em pointwise state constraints} written in the form
\begin{eqnarray}\label{state-cons}
\la x^j_*,x(t)\ra\le c_j\;\mbox{ for all }\;t\in[0,T]\;\mbox{ and }\;j=1,\ldots,s.
\end{eqnarray}\vspace*{-0.2in}

By a {\em feasible solution} to $(P)$ we understand a pair $(u(\cdot),x(\cdot))$ such that $u(\cdot)$ is measurable and that $x(\cdot)\in W^{1,2}([0,T],\R^n)$ subject to the constraints in \eqref{1.2}, \eqref{cont}, and hence in \eqref{state-cons}. Then \cite[Theorem~1]{et} implies that the set of feasible solutions to $(P)$ is nonempty under some assumptions that are much milder than those which are listed below.\vspace*{-0.05in}

Following \cite{cmn}, we say that a feasible pair $(\ox(\cdot),\ou(\cdot))$ for $(P)$ is a {\em $W^{1,2}\times L^2$-local minimizer} for this problem if there is $\e>0$ such that $J[\ox,\ou]\le J[x,u]$ for all the feasible pairs $(x(\cdot),u(\cdot))$ satisfying
\begin{equation*}
\int_0^T\(\n\dot{x}(t)-\dot{\ox}(t)\en^2+\n u(t)-\ou(t)\en^2\)dt<\e.
\end{equation*}
It is clear that this notion of local minimizers for $(P)$ includes, in the framework of sweeping control problems, {\em strong} ${\cal C}\times L^2$-local minimizers and occupies an intermediate position between the conventional notions of strong and weak minima in variational problems; cf.\  \cite{m06}.\vspace*{-0.05in}

Next we formulate the assumptions on the given data of $(P)$ needed for applications to the practical models considered below.
Note that the presented results taken from \cite{cmn} hold under more general assumptions, but we confine ourselves to the case of {\em smooth} functions and {\em convex} sets in $(P)$ that correspond to the models under consideration. In the following {\em standing assumptions} imposed in the rest of the paper without mentioning, the pair $(\ox(\cdot),\ou(\cdot))$ stands for the reference feasible solution to $(P)$, which is a chosen $W^{1,2}\times L^2$-local minimizer if stated so.\\[0.7ex]
{\bf(H1)} The control set $U$ is compact and convex in $\R^d$, and the image set $g(x,U)$ is convex in $\R^n$.\\
{\bf(H2)} The cost function $\ph\colon\R^n\to\R$ in \eqref{cost1} is ${\cal C}^1$-smooth around $\ox(T)$.\\
{\bf(H3)} The perturbation mapping $g\colon\R^n\times\R^d\to\R^n$ in \eqref{1.2} is ${\cal C}^1$-smooth around $(\ox(\cdot),\ou(\cdot))$ and satisfies the sublinear growth condition
\begin{equation*}
\|g(x,u)\|\le\be\big(1+\|x\|\big)\;\mbox{ for all }\;u\in U\;\mbox{ with some }\;\be>0.
\end{equation*}
{\bf(H4)} The vertices $x^j_*$ of \eqref{C} satisfy the linear independence constraint qualification
\begin{equation*}
\Big[\sum_{j\in I(\ox)}\al_j x^j_*=0,\;\al_j\in\R\Big]\Longrightarrow\big[\al_j=0\;\mbox{ for all }\;j\in I(\ox)\big\}
\end{equation*}
along the trajectory $\ox=\ox(t)$ as $t\in[0,T]$, where $I(\ox):=\{j\in\{1,\ldots,s\}\;|\;\la x^j_*,\ox\ra=c_j\}$.

First we present a crucial development of \cite{cmn} establishing close relationships between feasible and optimal solutions to problem $(P)$ and those to a sequence of its {\em discrete approximations}. Given any $m\in\N:=\{1,2,\ldots\}$, consider the discrete mesh
\begin{equation*}
\DD_m:=\big\{0=t_{0m}<t_{1m}<\ldots<t_{2^mm}=T\big\}\;\mbox{ with }\;h_m:=t_{(k+1)m}-t_{km}
\end{equation*}
on $[0,T]$ and the sequence of discrete-time inclusions approximating the controlled sweeping process \eqref{1.2}:
\begin{equation}\label{sw-disc}
x_{(k+1)m}\in x_{km}+h_m\big(g(x_{km},u_{km})-N(x_{km};C)\big)\;\mbox{ as }\;k=0,\ldots,2^m-1\;\mbox{ and }\;x_{0m}=x_0\in C
\end{equation}
over discrete pairs $(x_m,u_m)=(x_{0m},x_{1m},\ldots,x_{2^mm},u_{0m},u_{1m},\ldots,u_{(2^m-1)m})$ with the control constraints
\begin{equation}\label{cont-disc}
u_m=\big(u_{0m},u_{1m},\ldots,u_{(2^m-1)m}\big)\in U.
\end{equation}
Denote by $I_{km}:=[t_{(k-1)m},t_{km})$ for $k=1,\ldots,2^m$ the corresponding subintervals of $[0,T]$. The following theorem is a combination of the results taken from \cite[Theorems~3.1 and 4.2]{cmn}.\vspace*{-0.15in}

\begin{theorem}{\bf(discrete approximations in sweeping optimal control).}\label{Thm6.1*} Let $(\ox(\cdot),\ou(\cdot))$ be a feasible solution to problem $(P)$ such that $\ox(\cdot)\in W^{1,2}([0,T];\R^n)$ and that $\ou(\cdot)$ is of bounded variation $($BV$)$ with a right continuous representative on $[0,T]$. Then there exist sequences of unit vectors sequences $z_m^{jk}\to x^j_*$, vectors $c_m^{jk}\to c_j$ as $m\to\infty$, and state-control pairs $(\ox_m(t),\ou_m(t)),\;0\le t\le T$, for which we have:
{\bf (a)} The sequence of controls $\ou_m\colon[0,T]\to U$, which are constant on each interval $I_{km}$, converges to $\ou(\cdot)$
strongly in  $L^2([0,T];\R^d)$ and pointwise on $[0,T]$.\\
{\bf (b)} The sequence of continuous state mappings $\ox_m\colon[0,T]\to\R^n$, which are affine on each interval $I_{km}$, converges strongly in
$W^{1,2}([0,T];\R^n)$ to $\ox(\cdot)$, and satisfy the inclusions
\begin{equation*}
\ox_m(t_{km})=\ox(t_{km})\in C_{km}\;\mbox{ for each }\;k=1,\ldots,2^m\;\mbox{ with }\;\ox_m(0)=x_0,
\end{equation*}
where the perturbed polyhedra $C_{km}$ are given by
\begin{equation}\label{C-disc}
C_{km}:=\bigcap_{j=1}^s\big\{x\in\R^n\;\big|\;\langle z_m^{jk},x\rangle\le c_m^{jk}\big\}\;\mbox{ for }\;k=1,\ldots,2^m\;\mbox{ with }\;C_{0m}:=C.
\end{equation}
{\bf (c)} For all $t\in(t_{(k-1)m},t_{km})$ and $k=1,\ldots,2^m$ we have the differential inclusions
\begin{equation*}
\dot{\ox}_m(t)\in-N\big(\ox_m(t_{km});C_{km}\big)+g\big(\ox_m(t_{km}),\ou_m(t)\big).
\end{equation*}
If furthermore $(\ox(\cdot),\ou(\cdot))$ is a $W^{1,2}\times L^2$-local minimizer for problem $(P)$, then for each $m\in\N$ the pair $(\ox_m(\cdot),\ou_m(\cdot))$ above can be chosen so that its restriction on the discrete mesh $\DD_m$ is an optimal solution to the discrete sweeping control problem $(P_m)$ of minimizing the cost functional
\begin{equation*}
J_m[x_m,u_m]:=\vph\big(x_m(T)\big)+
\disp\frac{1}{2}\sum_{k=0}^{2^m-1}
\int_{t_{km}}^{t_{(k+1)m}}\Big(\Big\|\frac{x_{(k+1)m}-x_{km}}{h_m}-
\disp\dot{\ox}(t)\Big\|^2+\big\|u_{km}-\ou(t)\big\|^2\Big)dt
\end{equation*}
over all the pair $(x_m,u_m)$ satisfying \eqref{sw-disc}, \eqref{cont-disc}, $x_m(t_{km})\in C_{km}$ as $k=1,\ldots,2^m$ with $C_{km}$ taken from \eqref{C-disc}, and the $W^{1,2}\times L^2$-localization constraint
\begin{equation*}
\sum_{k=0}^{2^m-1}
\int_{t_{km}}^{t_{(k+1)m}}\Big(\Big\|\frac{x_{(k+1)m}-x_{km}}{h_m}-
\disp\dot{\ox}(t)\Big\|^2+\big\|u_{km}-\ou(t)\big\|^2\Big)dt\le\frac{\e}{2}.
\end{equation*}
\end{theorem}\vspace*{-0.05in}

Note that the results of \cite[Theorems~6.1 and 6.2]{cmn} contain necessary optimality conditions for the discrete control problems $(P_m)$ formulated in Theorem~\ref{Thm6.1*} that are not used in this paper. Nevertheless, they are very instrumental, together with the results of Theorem~\ref{Thm6.1*} above, to derive necessary optimality conditions for local minimizers of problem $(P)$, which are strongly employed in what follows. The next theorem presents these results in the case of the smoothness and convexity assumptions needed for the subsequent application to the practical models below; see \cite[Theorem~7.1]{cmn} for more general settings. Let us emphasize that, even in the case of smooth and convex data, the derivation of the obtained optimality conditions for $(P)$ is strongly based on the advanced tools of (nonconvex) first-order and second-order variational analysis and generalized differentiation taken from \cite{m-book}.\vspace*{-0.15in}

\begin{theorem}{\bf(necessary optimality conditions for controlled sweeping processes).}\label{Thm6.2*}
Let $(\ox(\cdot),\ou(\cdot))$ be a $W^{1,2}\times L^2$-local minimizer for $(P)$ under the assumptions of Theorem~{\rm\ref{Thm6.1*}}. Then there exist a multiplier $\lm\ge 0$, a measure $\gg=(\gg^1,\ldots,\gg^n)\in C^*([0,T];\R^n)$ as well as adjoint arcs $p(\cdot)\in W^{1,2}([0,T];\R^n)$ and $q(\cdot)\in BV([0,T];\R^n)$ such that $\lm+\|q(0\|+\|p(T)\|>0$ and the following conditions are satisfied:\\
$\bullet$ {\sc Primal velocity representation:}
\begin{eqnarray}\label{37}
-\dot{\ox}(t)=\sum_{j=1}^s\eta^j(t)x^j_*-g\big(\ox(t),\ou(t)\big)\;\mbox{ for a.e. }\;t\in[0,T],
\end{eqnarray}
where $\eta^j(\cdot)\in L^2([0,T];\R_+)$ being uniquely determined by  \eqref{37} and well defined at $t=T$.\\
$\bullet$ {\sc Adjoint system:}
\begin{eqnarray*}
\dot{p}(t)=-\nabla_x g\big(\ox(t),\ou(t)\big)^*q(t)\;\mbox{ for a.e. }\;t\in[0,T],
\end{eqnarray*}
where the dual arcs $q(\cdot)$ and $p(\cdot)$ are precisely connected by the equation
\begin{eqnarray*}
q(t)=p(t)-\int_{(t,T]}d\gg(\tau)
\end{eqnarray*}
that holds for all $t\in[0,T]$ except at most a countable subset.\\
$\bullet$ {\sc Maximization condition}:
\begin{eqnarray*}
\big\la\psi(t),\ou(t)\big\ra=\max\big\{\big\la\psi(t),u\big\ra\big|\;u\in U\big\}\;\mbox{ with }\;\psi(t):=\nabla_u g\big(\ox(t),\ou(t)\big)^*q(t)
\;\mbox{ for a.e. }\;t\in[0,T].
\end{eqnarray*}
$\bullet$ {\sc Complementarity conditions:}
\begin{eqnarray*}
\big\la x^j_*,\ox(t)\big\ra<c_j\Longrightarrow\eta^j(t)=0\;\mbox{ and }\;\eta^j(t)>0\Longrightarrow\big\la x^j_*,q(t)\big\ra=c_j
\end{eqnarray*}
for a.e.\ $t\in[0,T]$ including $t=T$ and for all $j=1,\ldots,s$.\\
$\bullet$ {\sc Right endpoint transversality conditions}:
\begin{eqnarray*}
-p(T)=\lm\nabla\vph\big(\ox(T)\big)+\sum_{j\in I(\ox(T))}\eta^j(T)x^j_*\;\mbox{ with }\;\sum_{j\in I(\ox(T))}\eta^j(T)x^j_*\in N\big(\ox(T);C\big).
\end{eqnarray*}
$\bullet$ {\sc Measure nonatomicity condition:} If $t\in[0,T)$ and $\la x^j_*,\ox(t)\ra<c_j$ for all $j=1,\ldots,s$, then there
is a neighborhood $V_t$ of $t$ in $[0,T]$ such that $\gg(V)=0$ for all the Borel subsets $V$ of $V_t$.
\end{theorem}\vspace*{-0.1in}

In the next two sections we develop applications of the obtained results to two classes of practical models formulated in the form of the sweeping optimal control problem $(P)$.\vspace*{-0.25in}

\section{Controlled Mobile Robot Model with Obstacles}\vspace*{-0.15in}

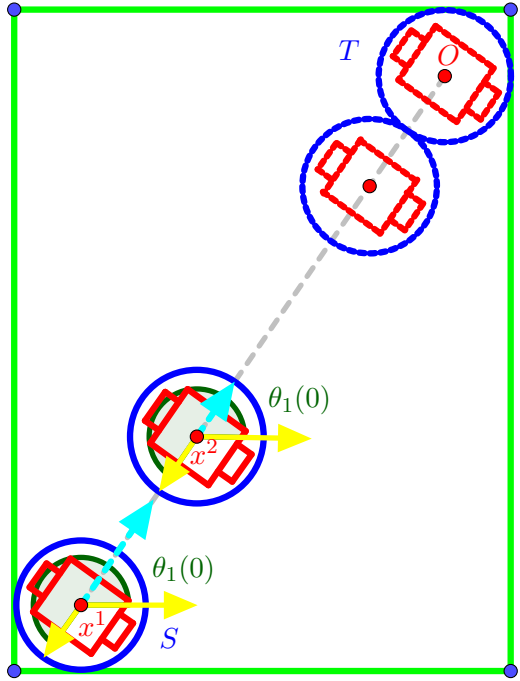
\begin{figure}[htp]
\begin{center}
\definecolor{qqwuqq}{rgb}{0.,0.39215686274509803,0.}
\definecolor{ffffqq}{rgb}{1.,1.,0.}
\definecolor{qqffff}{rgb}{0.,1.,1.}
\definecolor{cqcqcq}{rgb}{0.7529411764705882,0.7529411764705882,0.7529411764705882}
\definecolor{qqqqff}{rgb}{0.,0.,1.}
\definecolor{ffqqqq}{rgb}{1.,0.,0.}
\definecolor{qqffqq}{rgb}{0.,1.,0.}
\definecolor{ududff}{rgb}{0.30196078431372547,0.30196078431372547,1.}
\begin{tikzpicture}[line cap=round,line join=round,>=triangle 45,x=0.55cm,y=0.55cm]
\clip(-7.656,-9.041) rectangle (7.833,9.178);
\draw [shift={(-1.5852275239813505,-2.3298899439388006)},line width=2.pt,color=qqwuqq,fill=qqwuqq,fill opacity=0.10000000149011612] (0,0) -- (-1.0031717575446821:1.153636197971901) arc (-1.0031717575446821:235.50448252327277:1.153636197971901) -- cycle;
\draw [shift={(-4.385443499836288,-6.404912414635513)},line width=2.pt,color=qqwuqq,fill=qqwuqq,fill opacity=0.10000000149011612] (0,0) -- (-0.050775645832007464:1.153636197971901) arc (-0.050775645832007464:234.57272388924946:1.153636197971901) -- cycle;
\draw [line width=2.4pt,color=qqffqq] (-6.,8.)-- (6.,8.);
\draw [line width=2.4pt,color=qqffqq] (6.,8.)-- (6.,-8.);
\draw [line width=2.4pt,color=qqffqq] (6.,-8.)-- (-6.,-8.);
\draw [line width=2.4pt,color=qqffqq] (-6.,-8.)-- (-6.,8.);
\draw [line width=2.4pt,color=qqqqff] (-4.385443499836288,-6.404912414635513) circle (0.8598525918097203cm);
\draw [line width=2.4pt,dash pattern=on 2pt off 2pt,color=qqqqff] (4.407396921449587,6.390894268472447) circle (0.8791597246362263cm);
\draw [line width=2.pt,dash pattern=on 4pt off 4pt,color=cqcqcq] (4.407396921449587,6.390894268472447)-- (-4.385443499836288,-6.404912414635513);
\draw [line width=2.4pt,color=qqqqff] (-1.5852275239813505,-2.3298899439388006) circle (0.8885538363966832cm);
\draw [line width=2.4pt,dash pattern=on 2pt off 2pt,color=ffqqqq] (4.185071521922799,7.6034467516037845)-- (5.590795206225628,6.587000395262407);
\draw [line width=2.4pt,dash pattern=on 2pt off 2pt,color=ffqqqq] (3.286323521820033,6.297418771200554)-- (4.6920472061228615,5.280972414859177);
\draw [line width=2.4pt,color=ffqqqq] (-1.809089095175605,-1.1804831175311161)-- (-0.4033654108727758,-2.196929473872494);
\draw [line width=2.4pt,color=ffqqqq] (-2.712600595098011,-2.4889393887550697)-- (-1.3068769107951823,-3.5053857450964467);
\draw [line width=2.4pt,dash pattern=on 2pt off 2pt,color=ffqqqq] (4.185071521922799,7.6034467516037845)-- (3.286323521820033,6.297418771200554);
\draw [line width=2.4pt,dash pattern=on 2pt off 2pt,color=ffqqqq] (5.590795206225628,6.587000395262407)-- (4.6920472061228615,5.280972414859177);
\draw [line width=2.4pt,dash pattern=on 2pt off 2pt,color=ffqqqq] (5.3525049714216975,6.24072567300967)-- (5.718691429584372,5.992743738426331);
\draw [line width=2.4pt,dash pattern=on 2pt off 2pt,color=ffqqqq] (5.718691429584372,5.992743738426331)-- (5.274646177125915,5.341621734245304);
\draw [line width=2.4pt,dash pattern=on 2pt off 2pt,color=ffqqqq] (5.274646177125915,5.341621734245304)-- (4.911322674432222,5.599615568440938);
\draw [line width=2.4pt,dash pattern=on 2pt off 2pt,color=ffqqqq] (3.9434761539392293,7.2523691297495185)-- (3.5974594320077102,7.499300554886043);
\draw [line width=2.4pt,dash pattern=on 2pt off 2pt,color=ffqqqq] (3.5974594320077102,7.499300554886043)-- (3.1515186143353437,6.872572919238802);
\draw [line width=2.4pt,dash pattern=on 2pt off 2pt,color=ffqqqq] (3.1515186143353437,6.872572919238802)-- (3.5041404471896738,6.613942423109237);
\draw [line width=2.4pt,color=ffqqqq] (-4.641909518341075,-5.279206801376119)-- (-3.236185834038246,-6.295653157717496);
\draw [line width=2.4pt,color=ffqqqq] (-5.545329049997704,-6.568101754238313)-- (-4.139605365694875,-7.58454811057969);
\draw [line width=2.4pt,color=ffqqqq] (-2.712600595098011,-2.4889393887550697)-- (-1.809089095175605,-1.1804831175311161);
\draw [line width=2.4pt,color=ffqqqq] (-1.3068769107951823,-3.5053857450964467)-- (-0.4033654108727758,-2.196929473872494);
\draw [line width=2.4pt,color=ffqqqq] (-4.641909518341075,-5.279206801376119)-- (-5.545329049997704,-6.568101754238313);
\draw [line width=2.4pt,color=ffqqqq] (-3.236185834038246,-6.295653157717496)-- (-4.139605365694875,-7.58454811057969);
\draw [line width=2.4pt,color=ffqqqq] (-0.6175988640940816,-2.5071802071475595)-- (-0.2622910681792531,-2.7590883569354565);
\draw [line width=2.4pt,color=ffqqqq] (-0.2622910681792531,-2.7590883569354565)-- (-0.7303620177086542,-3.4189916628289314);
\draw [line width=2.4pt,color=ffqqqq] (-0.7303620177086542,-3.4189916628289314)-- (-1.0786673490731105,-3.1748949519918117);
\draw [line width=2.4pt,color=ffqqqq] (-2.020125069632086,-1.4861032942030847)-- (-2.3494271054251072,-1.2474493887841256);
\draw [line width=2.4pt,color=ffqqqq] (-2.3494271054251072,-1.2474493887841256)-- (-2.8405179377182495,-1.9457191659504767);
\draw [line width=2.4pt,color=ffqqqq] (-2.8405179377182495,-1.9457191659504767)-- (-2.4985200043398392,-2.1789100293623944);
\draw [line width=2.4pt,color=ffqqqq] (-3.459016949536382,-6.613562932791256)-- (-3.127045675697523,-6.843663397075664);
\draw [line width=2.4pt,color=ffqqqq] (-3.127045675697523,-6.843663397075664)-- (-3.5647617580901603,-7.45492561882335);
\draw [line width=2.4pt,color=ffqqqq] (-3.5647617580901603,-7.45492561882335)-- (-3.8887905947359784,-7.2267144784028305);
\draw [line width=2.4pt,color=ffqqqq] (-4.889141084847205,-5.631928330070142)-- (-5.218222190977455,-5.397608115639725);
\draw [line width=2.4pt,color=ffqqqq] (-5.218222190977455,-5.397608115639725)-- (-5.628638848454803,-5.975207864629086);
\draw [line width=2.4pt,color=ffqqqq] (-5.628638848454803,-5.975207864629086)-- (-5.296732755842964,-6.213433189292926);
\draw [->,line width=2.4pt,dash pattern=on 4pt off 4pt,color=qqffff] (-4.385443499836288,-6.404912414635513) -- (-2.629986263083249,-3.85027814792266);
\draw [->,line width=2.4pt,dash pattern=on 4pt off 4pt,color=qqffff] (-1.5852275239813505,-2.3298899439388006) -- (-0.6702727177479118,-0.9983992920668932);
\draw [->,line width=2.pt,color=ffffqq] (-1.5852275239813505,-2.3298899439388006) -- (1.190894310311549,-2.378501057252683);
\draw [->,line width=2.pt,color=ffffqq] (-4.385443499836288,-6.404912414635513) -- (-1.5673586038843745,-6.407409808322731);
\draw [->,line width=2.pt,color=ffffqq] (-1.5852275239813505,-2.3298899439388006) -- (-2.50018233021479,-3.6613805958107086);
\draw [->,line width=2.pt,color=ffffqq] (-4.385443499836288,-6.404912414635513) -- (-5.291679909109132,-7.678826133992425);
\draw [line width=2.4pt,dash pattern=on 2pt off 2pt,color=qqqqff] (2.5911674810734815,3.7287727457752515) circle (0.8942072796250358cm);
\draw [line width=2.4pt,dash pattern=on 2pt off 2pt,color=ffqqqq] (1.4607305252749676,3.592340699385865)-- (2.866454209577797,2.575894343044488);
\draw [line width=2.4pt,dash pattern=on 2pt off 2pt,color=ffqqqq] (2.3572839729772377,4.839719409231686)-- (3.7630076572800677,3.8232730528903094);
\draw [line width=2.4pt,dash pattern=on 2pt off 2pt,color=ffqqqq] (2.376774265318591,4.859209701573027)-- (1.478026265215825,3.553181721169796);
\draw [line width=2.4pt,dash pattern=on 2pt off 2pt,color=ffqqqq] (3.7410947292133487,3.8846950845059793)-- (2.8423467291105826,2.5786671041027485);
\draw [line width=2.4pt,dash pattern=on 2pt off 2pt,color=ffqqqq] (2.1097375195590855,4.5484364502892305)-- (1.7637207976275655,4.795367875425753);
\draw [line width=2.4pt,dash pattern=on 2pt off 2pt,color=ffqqqq] (1.7785682974959398,4.761001261742127)-- (1.3326274798235733,4.1342736260948865);
\draw [line width=2.4pt,dash pattern=on 2pt off 2pt,color=ffqqqq] (1.3391186118130278,4.134771249494191)-- (1.6917404446673598,3.8761407533646293);
\draw [line width=2.4pt,dash pattern=on 2pt off 2pt,color=ffqqqq] (3.885276540454481,3.3110142725813168)-- (3.441231287996024,2.6598922684002897);
\draw [line width=2.4pt,dash pattern=on 2pt off 2pt,color=ffqqqq] (3.423473386730296,2.6744747904213684)-- (3.0601498840366026,2.9324686246170026);
\draw [line width=2.4pt,dash pattern=on 2pt off 2pt,color=ffqqqq] (3.4896095549612967,3.571793702555571)-- (3.8557960131239706,3.3238117679722303);
\begin{large}
\draw [fill=ududff] (-6.,8.) circle (2.5pt);
\draw [fill=ududff] (6.,8.) circle (2.5pt);
\draw [fill=ududff] (6.,-8.) circle (2.5pt);
\draw [fill=ududff] (-6.,-8.) circle (2.5pt);
\draw [fill=ffqqqq] (-4.385443499836288,-6.404912414635513) circle (2.5pt);
\draw[color=ffqqqq] (-4.1,-6.85) node {$x^1$};
\draw[color=qqqqff] (-2.257,-7.2) node {$S$};
\draw [fill=ffqqqq] (4.407396921449587,6.390894268472447) circle (2.5pt);
\draw[color=ffqqqq] (4.5,6.9) node {$O$};
\draw[color=qqqqff] (2.1,7.03) node {$T$};
\draw [fill=ffqqqq] (-1.5852275239813505,-2.3298899439388006) circle (2.5pt);
\draw[color=ffqqqq] (-1.4,-2.7) node {$x^2$};
\draw[color=qqwuqq] (0.9,-1.45) node {$\th_1(0)$};
\draw[color=qqwuqq] (-1.9,-5.55) node {$\th_1(0)$};
\draw [fill=ffqqqq] (2.5911674810734815,3.7287727457752515) circle (2.5pt);
\end{large}
\end{tikzpicture}\caption{{\bf Mobile robot model with obstacles.}}
\label{pic1}
\end{center}
\end{figure}

In this section we formulate and investigate an {\em optimal control} version of the {\em mobile robot model} with obstacles which dynamics is described in \cite{HB} as a sweeping process. This model concerns $n$ mobile robots $(n\ge 2)$ identified with safety disks in the plane of the same radius $R$ as depicted in Fig.~\ref{pic1}.\vspace*{-0.05in}

The goal of each robot is to reach the target by the shortest path during a fixed time interval $[0,T]$ while avoiding the other $n-1$ robots that are treated by it as obstacles.\vspace*{-0.05in}

To formalize the model, consider the configuration vector $x=(x^1,\ldots,x^n)\in\R^{2n}$, where $x^i\in\R^2$ is the center of the safety disk $i$ with coordinates $(\|x^i\|\cos\th_i,\|x^i\|\sin\th_i)$. This means that the trajectory $x^i(t)$ of the $i$-robot/obstacle admits the representation
\begin{eqnarray*}
\ox^i(t)=\big(\|\ox^i(t)\|\cos\th_i(t),\|\ox^i(t)\|\sin\th_i(t)\big)\;\mbox{ for }\;i=1,\ldots,n,
\end{eqnarray*}
where the angle $\th_i$ signifies the corresponding direction. According to the model dynamics, at the moment of contacting the obstacle (one or more) the robot in question keeps its velocity and pushes the other robots in contact to go to the target with the same velocity and then to maintain their constant velocities until reaching either other obstacles or the end of the process at the final time $t=T$. In this framework, the constant direction $\th_i$ of $x^i$ is the smallest positive angle in standard position formed by the positive $x$-axis and $Ox^i$; see Fig.~\ref{pic1}, where the origin is the target point.

To ensure the avoidance of collision between the robot and obstacles, we define the {\em admissible configuration} set by imposing the
{\em noncollision/nonoverlapping conditions} $\|x^i-x^j\|\ge 2R$ formulated as
\begin{eqnarray}\label{t:96}
Q_0:=\big\{x=\(x^1,\ldots,x^n\)\in\mathbb{R}^{2n}\big|\;D_{ij}(x)\ge 0\;\mbox{ whenever }\;i,j\in\{1,\ldots,n\}\big\},
\end{eqnarray}
where $D_{ij}(x)=\|x^{i}-x^j\|-2R$ is the distance between the safety disks $i$ and $j$.\vspace*{-0.05in}

Let $\nabla D_{ij}(x)$ be the gradient of $D_{ij}(x)$ at $x\ne 0$. In order to efficiently describe nonoverlapping of the safety disks, define the set of {\em admissible velocities} by
\begin{eqnarray*}
V_h(x):=\big\{v\in\R^{2n}\big|\;D_{ij}(x)+h\nabla D_{ij}(x)v\ge 0\;\mbox{ for all }\;i,j\in\{1,\ldots,n\},\;i<j\big\},\quad x\in\R^{2n},
\end{eqnarray*}
which is closely related to the admissible configuration set \eqref{t:96}. Indeed, if the chosen admissible configuration at time $t_k\in[0,T]$ is $x_k:=x(t_k)\in Q_0$, then the next configuration after the period of time $h>0$ is $x_{k+1}=x(t_k+h)$. Thus it follows from the first-order Taylor expansion at $x_k\ne 0$ that
\begin{eqnarray}\label{h:2.3}
D_{ij}\big(x(t_k+h)\big)=D_{ij}\big(x(t_k)\big)+h\nabla D_{ij}\big(x(t_k)\big)\dot{x}(t_k)+o(h)\;\mbox{ for small }\;h>0.
\end{eqnarray}
Taking now the admissible velocity $\dot{x}(t_k)\in V_h(x_k)$ and ignoring the term $o(h)$ for small $h$ give us
\begin{eqnarray*}
D_{ij}(x_k)+h\big\la\nabla D_{ij}(x_k),\dot{x}(t_k)\big\ra\ge 0,
\end{eqnarray*}
and therefore it follows from \eqref{h:2.3} that $D_{ij}(x(t_k+h))\ge 0$, i.e., $x(t_k+h)\in Q_0$.\vspace*{-0.05in}

Since all the robots intend to reach the target by the {\em shortest path}, their desired spontaneous (i.e., in the absence of other robots) velocities can be represented as
\begin{equation*}
S(x)=\big(S_0(x^1),\ldots,S_0(x^n)\big)\;\mbox{ with }\;S_0(x)=-s_0D(x),
\end{equation*}
where $D(x)$ stands for the distance from the position $x=(x^1,\ldots,x^n)\in Q_0$ to the target, and where the scalar $s_0\ge0$ indicates the speed. Due to $x\ne 0$ and hence by $\|D(x)\|\ne 1$, we get $s_0=\|S_0(x)\|$. Remembering that in the absence of obstacles the robots tend to keep their desired spontaneous velocities till reaching the target and taking into account the previous discussions, we describe the velocities by
\begin{equation*}
g\big(x(t)\big):=-\big(s_1\cos\th_1,s_1\sin\th_1,\ldots,s_n\cos\th_n,s_n\sin\th_n\big)\in\R^{2n}
\end{equation*}
for all $x\in Q_0$, where $s_i$ denotes the speed of robot $i$. However, if the robot in question touches the obstacles in the sense that $\|x^{i}(t)-x^1(t)\|=2R$, its velocity should be adjusted in order to keep the distance to be at least $2R$ by using some {\em control actions} in the velocity term. It can be modeled as
\begin{eqnarray}\label{t:99**}
g\big(x(t),u(t)\big)=\big(s_1u^1(t)\cos\th_1(t),s_1u^1(t)\sin\th_1(t),\ldots,s_nu^n(t)\cos\th_n(t),s_nu^n(t)\sin\th_n(t)\big)
\end{eqnarray}
with practically motivated control constraints represented by
\begin{eqnarray}\label{t:99*}
u(t)=\big(u^1(t),\ldots,u^n(t)\big)\in U\;\mbox{ for a.e. }t\in[0,T],
\end{eqnarray}
where the control set $U\subset\R^n$ will be specified below in particular settings.\vspace*{-0.05in}

To avoid overlapping between the robot in question and obstacles, we proceed as follows. Taking $x_k\in Q_0$ as the admissible configuration at the time $t_k$ and using the mapping $g\colon\R^{2n}\times\R^n\to\R^{2n}$ from \eqref{t:99**} with a given feasible control $u_k:=u(t_k)$ from \eqref{t:99*}, the next configuration $x_{k+1}$ is calculated by
\begin{eqnarray}\label{design}
x_{k+1}=x_k+h V_{k+1},
\end{eqnarray}
where $V_{k+1}\in\R^{2n}$ solves the {\em convex optimization problem}:
\begin{eqnarray}\label{Pi}
\mbox{minimize }\;\|V-g(x_k,u_k)\|^2\;\mbox{ subject to }\;V\in V_{h}(x_k),
\end{eqnarray}
and where the control $u_k\in U$ is involved into the desired velocity term to adjust the actual velocities of the robots and make sure that they do not overlap. The algorithmic design in \eqref{design} and \eqref{Pi} means therefore that $V_{k+1}$ is selected as the (unique) element from the set of admissible velocities as the closest one to the desired velocity $g(x_k,u_k)$ in order to avoid the robot overlapping.\vspace*{-0.05in}

Fix further any $m\in\N$ and divide $[0,T]$ into the $2^m$ equal subintervals of length $h_m:=T/2^m\dn 0$ as $m\to\infty$. Invoking the discrete time $t_{km}:=kh_m$, denote $I_{km}:=[t_{km},t_{(k+1)m})$ for $k=0,\ldots,2^m-1$ and $I_{2^m m}:=\{T\}$. Then according to \eqref{design} and \eqref{Pi} we have the algorithm
\begin{eqnarray}\label{alg}
x_{0m}\in Q_0\;\mbox{ and }\;x_{(k+1)m}:=x_{km}+h_m V_{(k+1)m}\;\mbox{ for all }\;k=0,\ldots,2^m-1,
\end{eqnarray}
where $V_{(k+1)m}$ is defined as the {\em projection} of $g(x_{km},u_{km})$ onto the admissible velocity set $V_{h_m}(x_{km})$ by
\begin{eqnarray}\label{vk}
V_{(k+1)m}:=\disp\Pi\big(g(x_{km},u_{km});V_{h_m}(x_{km})\big),\quad k=0,\ldots,2^m-1.
\end{eqnarray}
Invoking the construction of $x_{km}$ for $0\le k\le 2^m-1$ and $m\in\N$, define next a sequence of piecewise linear mappings $x_{2^m}\colon[0,T]\to\R^{2n}$, $m\in\N$, which pass through those points by:
\begin{eqnarray}\label{xk}
x_{2^m}(t):=x_{km}+(t-t_{km})V_{(k+1)m}\;\mbox{ for all }\;t\in I_{km},\quad k=0,\ldots,2^m-1.
\end{eqnarray}
Whenever $m\in\N$, we clearly have the relationships
\begin{eqnarray}\label{2m}
x_{2^m}(t_{km})=x_{km}=\underset{t\to t_{km}}{\lim}x_{km}(t)\;\mbox{ and }\;\dot{x}_{2^m}(t):=V_{(k+1)m}\;\mbox{ for all }\;t\in(t_{km},t_{(k+1)m}).
\end{eqnarray}\vspace*{-0.2in}

As discussed in \cite{HB}, based on the results of \cite{v}, the solutions to \eqref{xk} in the {\em uncontrolled} setting of \eqref{vk} with $g=g(x)$ uniformly converge on $[0,T]$ to a trajectory of a certain perturbed sweeping process. The {\em controlled} model under consideration here is significantly more involved. In order to proceed by using the results of Theorem~\ref{Thm6.1*}, for all $x\in\R^{2n}$ consider the set
\begin{eqnarray}\label{K}
K(x):=\big\{y\in\R^{2n}\big|\;D_{ij}(x)+\nabla D_{ij}(x)(y-x)\ge 0\;\mbox{ whenever }\;i<j\big\},
\end{eqnarray}
which allows us to represent the algorithm in \eqref{vk}, \eqref{xk} as
\begin{eqnarray*}
x_{(k+1)m}=\Pi\Big(x_{km}+h_m g(x_{km},u_{km});K(x_{km})\Big)\;\mbox{ for }\;k=0,\ldots,2^m-1.
\end{eqnarray*}
It can be equivalently rewritten in the form
\begin{eqnarray*}
x_{2^m}\big(\vt_{2^m}(t)\big)=\Pi\Big(x_{2^m}(\tau_{2^m}(t))+h_m g\big(x_{2^m}(\tau_{2^m}(t)),u_{2^m}(\tau_{2^m}(t)\big);K(x_{2^m}(\tau_{2^m}(t))\Big)\;\mbox{ for all }\;t\in[0,T],
\end{eqnarray*}
where the functions $\tau_{2^m}(\cdot)$ and $\vt_{2^m}(\cdot)$ are defined by $\tau_{2^m}(t):=t_{km}$ and $\vt_{2^m}(t):=t_{(k+1)m}$ for all $t\in I_{km}$. Taking into account the construction of the convex set $K(x)$ in \eqref{K} and definition \eqref{1.4} of the normal cone together with the relationships in \eqref{2m}, we arrive at the sweeping process inclusions
\begin{eqnarray}\label{h:3.2}
\dot x_{2^m}(t)\in-N\big(x_{2^m}(\vt_{2^m}(t));K(x_{2^m}(\tau_{2^m}(t)))\big)+g\big(x_{2^m}(\tau_{2^m}(t)),u_{2^m}(\tau_{2^m}(t))\big)\;\mbox{ a.e. }\;t\in[0,T]
\end{eqnarray}
with $x_{2^m}(0)=x_0\in K(x_0)=Q_0$ and $x_{2^m}(\vt_{2^m}(t))\in K(x_{2^m}(\tau_{2^m}(t)))$ on $[0,T]$. To formalize \eqref{h:3.2} as a controlled perturbed sweeping process of type \eqref{1.2}, define the convex polyhedron $C\subset\R^{2n}$ as in \eqref{C} by
\begin{eqnarray}\label{e:131***}
C:=\big\{x\in\R^{2n}\big|\;\la x^j_*,x\ra\le c_j,\;j=1,\ldots,n-1\big\}
\end{eqnarray}
with $c_j:=-2R$ and with the $n-1$ vertices of the polyhedron given by
\begin{eqnarray}\label{e}
x^j_*:=e_{j1}+e_{j2}-e_{(j+1)1}-e_{(j+1)2}\in\R^{2n},\quad j=1,\ldots,n-1,
\end{eqnarray}
where $e_{ji}$ for $j=1,\ldots,n$ and $i=1,2$ are the vectors in $\R^{2n}$ of the form
\begin{eqnarray*}
e:=\big(e_{11},e_{12},e_{21},e_{22},\ldots,e_{n1},e_{n2}\big)\in\R^{2n}
\end{eqnarray*}
with 1 at only one position of $e_{ji}$ and $0$ at all the other positions.\vspace*{-0.05in}

We now formulate the {\em sweeping optimal control problem} of type $(P)$ from Section~2 that can be treated as a continuous-time counterpart of the discrete algorithm of the controlled mobile robot model by taking into account the model goal stated above. Consider the cost functional
\begin{equation}\label{t:102*}
\mbox{minimize }\;J[x,u]:=\disp\frac{1}{2}\big\|x(T)\big\|^2,
\end{equation}
which reflects model goal to {\em minimize the distance} of the robot from the admissible configuration set to the target. We describe the continuous-time dynamics by the controlled sweeping process
\begin{eqnarray}\label{t:101*}
\left\{\begin{array}{lcl}
-\dot{x}(t)\in N\big(x(t);C\big)-g\big(x(t),u(t)\big)\;\mbox{ for a.e. }\;t\in[0,T],\\
x(0)=x_0\in C,\;u(t)\in U\;\mbox{ a.e. on }[0,T],
\end{array}\right.
\end{eqnarray}
where the constant set $C$ is taken from \eqref{e:131***}, the control constraints reduce to \eqref{t:99*}, and the dynamic noncollision condition $\|x^i(t)-x^j(t)\|\ge 2R$ amounts to the pointwise state constraints
\begin{eqnarray}\label{t:100*}
x(t)\in C\Longleftrightarrow\la x^j_*,x(t)\ra\le c_j\;\mbox{ for all }\;t\in[0,T]\;\mbox{ and }\;j=1,\ldots,n-1,
\end{eqnarray}
which follow from \eqref{t:101*} due to the construction of $C$ and the normal cone definition \eqref{1.4}. Next we obtain other representations of $C$, which allow us to make connections between the discrete dynamics in \eqref{vk}, \eqref{xk} and its sweeping control counterpart in \eqref{t:101*}. Taking into account that we are interested in the limiting process when the discrete step $h$ in \eqref{design} diminishes, it is possible to choose in what follows a convenient equivalent norm $\|(x^{j1},x^{j2})\|:=|x^{j1}|+|x^{j2}|$ for each component $x^j\in\R^2$ of $x\in\R^{2n}$.\vspace*{-0.15in}

\begin{lemma}{\bf(sweeping set representations).}\label{sw-rep} In addition to the noncollision conditions
\begin{equation}\label{noncol-disc}
\|x^i_{km}-x^j_{km}\|\ge 2R\;\mbox{ for all }\;i,j\in\big\{1,\ldots,n\big\},\;k=0,\ldots,2^m-1,\;\mbox{ and }\;m\in\N
\end{equation}
imposed on the points $x_{km}$ from \eqref{t:96}, suppose that
\begin{equation}\label{comp-disc}
x^{(j+1)1}_{km}>x^{j1}_{km}\;\mbox{ and }\;x^{(j+1)2}_{km}>x^{j2}_{km}\;\mbox{ whenever }\;j=1,\ldots,n-1
\end{equation}
for the components $(x^{j1}_{km},x^{j2}_{km})\in\R^2$ of the iterations above as $k=0,\ldots,2^m-1$ and $m\in\N$.
Then the sweeping set $C$ in \eqref{e:131***} admits the representations
\begin{equation}\label{sw-rep1}
C=Q_0=K(x_{km})\;\mbox{ whenever }\;\;k=0,\ldots,2^m-1\;\mbox{ and }\;m\in\N,
\end{equation}
where the sets $Q_0$ and $K(x_{km})$ are defined by \eqref{t:96} and \eqref{K}, respectively.
\end{lemma}\vspace*{-0.15in}
{\bf Proof.} By using the noncollision conditions \eqref{noncol-disc} as well as the component conditions imposed in \eqref{comp-disc}, we can easily verify the relationships
\begin{equation*}
\begin{array}{ll}
\big\|x^j-x^{j+1}_{km}\big\|=\big|x^{j1}_{km}-x^{(j+1)1}_{km}\big|+\big|x^{j2}_{km}-x^{(j+1)2}_{km}\big|\\
=-x^{j1}_{km}-x^{j2}_{km}+x^{(j+1)1}_{km}+x^{(j+1)2}_{km}\ge 2R,
\end{array}
\end{equation*}
which yield the following equalities, where we use \eqref{e} and put $x:=x_{km}$ for the simplicity of notation:
\begin{equation*}
\begin{array}{ll}
C&=\big\{x\in\R^{2n}\big|\;\la x^j_*,x\ra\le c_j,\;j=1,\ldots,n-1\big\}\\
&=\big\{x\in\R^{2n}\big|\;x^{j1}+x^{j2}-x^{(j+1)1}-x^{(j+1)2}\le-2R,\;j=1,\ldots,n-1\big\}\\
&=\big\{x\in\R^{2n}\big|\;-x^{j1}-x^{j2}+x^{(j+1)1}+x^{(j+1)2}\ge 2R,\;j=1,\ldots,n-1\big\}\\
&=\big\{x\in\R^{2n}\big|\;\|x^{j+1}-x^j\|\ge 2R,\;j=1,\ldots,n-1\big\}.
\end{array}
\end{equation*}
Then it follows from the definition of $Q_0$ in \eqref{t:96} and $D_{ij}$ therein that $C=Q_0$.\vspace*{-0.05in}

To verify the second equality(ies) in \eqref{sw-rep1} for all the indicated indices $j$ and $m$ therein, we get by constructions of $K(\cdot)$ and $D_{ij}(\cdot)$ that
\begin{equation*}
\begin{array}{ll}
K(x_0)&=\big\{x\in\R^{2n}\big|\;D_{ij}(x_{0})+\nabla D_{ij}(x_{0})(x-x_{0})\ge 0\;\mbox{ if }\;i<j\big\},\\
&=\big\{x\in\R^{2n}\big|\;D_{ij}(x)\ge 0\;\mbox{ if }\;i<j\big\}=C,
\end{array}
\end{equation*}
where we drop indicating the dependence on the vector $x:=x_{km}$ from the indices $km$ as above. Thus we directly arrive at the second statement in \eqref{sw-rep1} and complete the proof of the lemma. $\h$\vspace*{-0.02in}

It follows from the defined constructions that we can replace $K(x_{2^m}(\tau_{2^m}(t)))$ by $C$ on $[0,T]$ for large $m$. Thus the sweeping process in \eqref{t:101*} can be treated as the limiting case of \eqref{h:3.2}. The next theorem provides an application and a specification of Theorem~\ref{Thm6.1*} for the robotics model under consideration.\vspace*{-0.12in}

\begin{theorem}{\bf(sweeping process description of the controlled mobile robot model).}\label{sw-robot} Let the pair $(\ox(\cdot),\ou(\cdot))$ satisfy the controlled sweeping system \eqref{t:101*}, where $C$ is taken from \eqref{e:131***}, $g$ is defined in \eqref{t:99**}, $U\subset\R^n$ is compact and convex, and the conditions
\begin{equation}\label{comp-cont}
\ox^{(j+1)1}(t)>\ox^{j1}(t)\;\mbox{ and }\;\ox^{(j+1)2}(t)>\ox^{j2}(t)\;\mbox{ for all }\;j=1,\ldots,n-1\;\mbox{ and }\;t\in[0,T]
\end{equation}
are fulfilled. Assume that $\ox(\cdot)\in W^{1,2}([0,T];\R^{2n})$ and that $\ou(\cdot)$ is BV on $[0,T]$ with a right continuous representative. Then there exist a sequence of state-control pairs $(\ox_m(t),\ou_m(t)),\;0\le t\le T$, satisfying \eqref{h:3.2} with $K(x_{2^m}(\tau_{2^m}(t))\equiv C$ for which all the conclusions of Theorem~{\rm\ref{Thm6.1*}} hold with $s=n-1$ and $C_{km}\equiv C$. If furthermore $(\ox(\cdot),\ou(\cdot))$ is a $W^{1,2}\times L^2$-local minimizer of the cost functional \eqref{t:102*} over the constrained dynamics \eqref{t:101*}, then any sequence of the $($extended on $[0,T])$ optimal solutions $(\ox_m(\cdot),\ou_m(\cdot))$ to the corresponding specifications of problems $(P_m)$ from Theorem~{\rm\ref{Thm6.1*}} converges to $(\ox(\cdot),\ou(\cdot))$ in the norm topology of $W^{1,2}([0,T];\R^{2n})\times L^2([0,T];\R^n)$.
\end{theorem}\vspace*{-0.12in}
{\bf Proof.} It can be directly checked that the assumptions (H1)--(H4) are satisfied in the setting of \eqref{t:102*}, \eqref{t:101*} with the data specified for the mobile robot model. Then using Lemma~\ref{sw-rep} and invoking the above discussions, we deduce the conclusions of the theorem from the corresponding results of Theorem~{\rm\ref{Thm6.1*}}. $\h$\vspace*{-0.02in}

From now on in this section, we exclusively study the continuous-time sweeping optimal control problem defined in \eqref{t:102*} and \eqref{t:101*} with the mobile robot model data. We label this problem as $(SR)$. Applying Theorem~\ref{Thm6.2*} allows us to obtain the following necessary optimality conditions for problem $(SR)$ that are formulated entirely in terms of the model data.\vspace*{-0.15in}

\begin{theorem}{\bf(necessary optimality conditions for the sweeping controlled robot model).}\label{Thm6.2*a} Let $(\ox(\cdot),\ou(\cdot))$ be a $W^{1,2}\times L^2$-local minimizer for problem $(SR)$, and let all the assumptions of Theorem~{\rm\ref{sw-robot}} be fulfilled.
Then there exist a multiplier $\lm\ge 0$, a measure $\gg\in C^*([0,T];\R^{2n})$ as well as adjoint arcs $p(\cdot)\in W^{1,2}([0,T];\R^{2n})$ and $q(\cdot)\in BV([0,T];\R^{2n})$ satisfying to the conditions:\vspace*{-0.1in}
\begin{itemize}
\item[\bf(1)] $-\dot{\ox}(t)=\disp{\sum^{n-1}_{j=1}}\eta^j(t)x_*^j-\big(g(\ox^1(t),\ou^1(t))),\ldots,g(\ox^n(t),\ou^n(t)\big)$ for
a.e.\ $t\in[0,T]$,\\ where $\eta^j(\cdot)\in L^2([0,T];\R_+)$ are uniquely defined by this representation and well defined at $t=T$;
\item[\bf(2)] $\|\ox^{j}(t)-\ox^{j+1}(t)\|>2R\Longrightarrow\eta^j(t)=0$ for all $j=1,\ldots,n-1$ and a.e.\ $t\in[0,T]$ including $t=T$;
\item[\bf(3)] $\eta^j(t)>0\Longrightarrow\la x^j_*,q(t)\ra=c_j$ for all $j=1,\ldots,n-1$ and a.e.\ $t\in[0,T]$ including $t=T$;
\item[\bf(4)] $\dot{p}(t)=-\nabla_x g\big(\ox(t),\ou(t)\big)^*q(t)\;\mbox{ for a.e. }t\in[0,T]$;
\item[\bf(5)] $q(t)=p(t)-\gamma([t,T])$ for all $t\in[0,T]$ except at most a countable subset;
\item[\bf(6)] $\big\la\psi(t),\ou(t)\big\ra=\max_{u\in U}\big\la\psi(t),u\big\ra\;\mbox{ for a.e. }\;t\in[0,T]$, where $\psi(t):=
\nabla_u g\big(\ox(t),\ou(t)\big)^*q(t)$;
\item[\bf(7)] $-p(T)=\lm\ox(T)+\sum_{j\in I(\ox(T))}\eta^j(T)x^j_*$ via the set of active constraint indices $I(\ox(T))$ at $\ox(T)$;
\item[\bf(8)] $\sum_{j\in I(\ox(T))}\eta^j(T)x^j_*\in N\big(\ox(T);C)$;
\item[\bf(9)] $\lm+\|q(0)\|+\|p(T)\|>0$.
\end{itemize}
\end{theorem}\vspace*{-0.15in}
{\bf Proof.} As discussed in the framework of Theorem~\ref{Thm6.2*a}, all the assumptions of  Theorem~\ref{Thm6.2*} are fulfilled for $(SR)$. Thus we can apply to the given local minimizer $(\ox(\cdot),\ou(\cdot))$ of $(SR)$ the necessary optimality conditions from that theorem, which are specified as (1)--(9) in the setting under consideration. The only thing we need to check is the validity of the implication in (2). Indeed, we have the implication $\la\ox(t),x^j_*\ra<c_j\Longrightarrow\eta^j(t)=0$ for $j=1,\ldots,n-1$ and a.e. $t\in[0,T]$. Moreover, the conditions $\la\ox(t),x^j_*\ra<c_j$ are equivalent to
\begin{equation*}
-\ox^{j1}(t)-\ox^{j2}(t)+\ox^{(j+1)1}(t)+\ox^{(j+1)2}(t)>2R\;\mbox{ whenever }\;j=1,\ldots,n-1,\;\mbox{ a.e. on }\;[0,T].
\end{equation*}
By \eqref{comp-cont} the latter conditions are equivalent in turn to those in $\|\ox^j(t)-\ox^{j+1}(t)\|>2R$ for all $j,t$ indicated above. To verify this, we get while remembering the sum norm under consideration that
\begin{equation*}
\|\ox^j(t)-\ox^{j+1}(t)\|=|\ox^{j}(t)-\ox^{(j+1)1}(t)|+|\ox^{j2}(t)-\ox^{(j+1)2}(t)|.
\end{equation*}
Finally, it allows us to obtain the relationships
\begin{equation*}
\begin{array}{ll}
\|\ox^j(t)-\ox^{j+1}(t)\|=|\ox^{j1}(t)-\ox^{(j+1)1}(t)|+\big|\ox^{j2}(t)-\ox^{(j+1)2}(t)|\\
=-\ox^{j1}(t)-\ox^{j2}(t)+\ox^{(j+1)1}(t)+\ox^{(j+1)2}(t)>2R,
\end{array}
\end{equation*}
which justify (2) and thus completes the proof of the theorem. $\h$

Let us now discuss some conclusions for the mobile robot model that can be derived from the obtained theorem by taking into account the specific form of the perturbation mapping $g$ in \eqref{t:99**}.\vspace*{-0.05in}

$\bullet$ We know from the model description that at the contact time $t_1\in[0,T]$ when $\|\ox^{i}(t_1)-\ox^1(t_1)\|=2R$ for some $i=2,\ldots,n$, the robot in question tends to adjust its velocity in order to keep the distance between the obstacle in contact to be at least $2R$. By the model requirement, the robot maintains its constant velocity after the time $t=t_1$ until either reaching other obstacles ahead or stopping at $t=T$. If furthermore the robot touches other obstacles, it pushes them to go to the target in the same direction as before $t=t_1$. By \eqref{t:99**}, the differential relation in (1) is written for a.e.\ $t\in[0,T]$ as
\begin{eqnarray}\label{t:103*}
\left\{\begin{array}{ll}
-\big(\dot{\ox}^{11}(t),\dot{\ox}^{12}(t)\big)=\eta^1(t)(1,1)-\big(s_1\ou^1(t)\cos\th_1(t),s_1\ou^1(t)\sin\th_1(t)\big),\\
-\big(\dot{\ox}^{i1}(t),\dot{\ox}^{i2}(t)\big)=\eta^{i-1}(t)(-1,-1)+\eta^{i}(t)(1,1)-\big(s_i\ou^i(t)\cos\th_i(t),s_i\ou^i(t)\sin\th_i(t)\big)\\
-(\dot{\ox}^{n1}(t),\dot{\ox}^{n2}(t)\big)=\eta^{n-1}(t)(-1,-1)-\big(s_n\ou^n(t)\cos\th_n(t),s_n\ou^n(t)
\sin\th_n(t)\big).
\end{array}\right.
\end{eqnarray}\vspace*{-0.1in}

$\bullet$ If the robot under consideration (robot~1) does not touch the first obstacle (robot~2) in the sense that $\|\ox^2(t)-\ox^1(t)\|>2R$ for all $t\in[0,T]$, then we get from (2) of Theorem~\ref{Thm6.2*a} that
\begin{eqnarray*}
\|\ox^{2}(t)-\ox^{1}(t)\|>2R\Longrightarrow\eta^1(t)=0\;\mbox{ for a.e. }\;t\in[0,T]\;\mbox{ including }\;t=T.
\end{eqnarray*}
Plugging $\eta^1(t)=0$ into \eqref{t:103*} gives us the equation
\begin{eqnarray*}
-\big(\dot{\ox}^{11}(t),\dot{\ox}^{12}(t)\big)=-\big(s_1\ou^1(t)\cos\th_1(t),s_1\ou^1(t)\sin\th_1(t)\big)\;\mbox{ a.e. on }\;[0,T],
\end{eqnarray*}
which means that the actual velocity and the spontaneous velocity of the robot agree for a.e.\ $t\in[0,T]$. Similarly we conclude that the condition $\|\ox^n(t)-\ox^{n-1}(t)\|>2R$ on $[0,T]$ yields $-\dot{\ox}^n(t)=-g(\ox^n(t),\ou^n(t))$ for a.e.\ $t\in[0,T]$, and then continue in this way with robots $i$.\vspace*{-0.05in}

To proceed further, assume that $\lm=1$ (otherwise we do not have enough information to efficiently employ Theorem~\ref{Thm6.2*a}) and also suppose for simplicity to handle the examples below that the control actions $\ou^i(\cdot)$ are constant on $[0,T]$ for all $i=1,\ldots,n$. Applying the Newton-Leibniz formula in \eqref{t:103*} gives us the trajectory representations
\begin{eqnarray}\label{nlf}
\left\{\begin{array}{ll}
\(\ox^{11}(t),\ox^{12}(t)\)=\(x^{11}_0,x^{12}_0\)-\disp\int^t_0\eta^1(\tau)\(1,1\)d\tau+t\(s_1\ou^1\cos\th_1,s_1\ou^1\sin\th_1\),\\
\(\ox^{i1}(t),\ox^{i2}(t)\)=\(x^{i1}_0,x^{i2}_0\)+\disp\int^t_0\eta^{i-1}(\tau)\(1,1\)d\tau-\disp\int^t_0\eta^{i}(\tau)\(1,1\)d\tau\\
\quad\quad\quad\quad\quad\quad\quad+t\(s_i\ou^i\cos\th_i,s_i\ou^i\sin\th_i\)\;\mbox{ whenever }\;i=2,\ldots,n-1,\\
\(\ox^{n1}(t),\ox^{n2}(t)\)=\(x^{n1}_0,x^{n2}_0\)+\disp\int^t_0\eta^{n-1}(\tau)\(1,1\)d\tau+t\(s_n\ou^n\cos\th_n,s_n\ou^n\sin\th_n\)
\end{array}\right.
\end{eqnarray}
for all $t\in[0,T]$, where $x_0:=(x^{11}_0,x^{12}_0\ldots,x^{n1}_0,x^{n2}_0)\in C$ stands for the starting point in \eqref{t:101*}.\vspace*{-0.05in}

Next we employ the obtained necessary optimality conditions to find optimal solutions and understand the sweeping process behavior in some typical situations that appear in the controlled robot mobile model by considering for simplicity the case of $n=2$. Note that in all the cases below we have the existence of optimal solutions by \cite{et}, and thus the unique ones determined by using necessary optimality conditions are indeed globally optimal for this model in the settings under consideration.\vspace*{-0.05in}

$\textbf{A: Mobile robot model without changing direction in contact.}$ The first typical situation is when the robot in question touches the other robot
(obstacle) so that there is no change of direction at the point of contact; see Fig.~\ref{pic1}. Let $t_1$ be the contact time, i.e.,
\begin{eqnarray}\label{t-contact}
t_1:=\min\big\{t\in[0,T]\big|\;\;\|\ox^1(t)-\ox^2(t)\|=2R\big\}.
\end{eqnarray}
Recalling that in our model the equal angles $\th_1=\th_2$ of the robot directions are constant together with the optimal controls, we get by \eqref{t:103*} the dynamic equations prior to and after time $t_1$:
\begin{eqnarray*}
\left\{\begin{array}{ll}
\dot{\ox}^1(t)=\(s_1\ou^1\cos\th_1,s_1\ou^1\sin\th_1\),\\
\dot{\ox}^2(t)=\(s_2\ou^2\cos\th_1,s_2\ou^2\sin\th_1)\)\;\mbox{ for }\;t\in[0,t_1),
\end{array}\right.
\end{eqnarray*}\vspace*{-0.2in}
\begin{eqnarray*}
\left\{\begin{array}{ll}
\dot{\ox}^1(t)=\(-\eta^1(t)+s_1\ou^1(t)\cos\th_1,-\eta^1(t)+s_1\ou^1\sin\th_1\),\\
\dot{\ox}^2(t)=\(\eta^1(t)+s_2\ou^2\cos\th_1,\eta^1(t)+s_2\ou^2\sin\th_1\)\;\mbox{ for }\;t\in[t_1,T].
\end{array}\right.
\end{eqnarray*}
This implies, with taking into account condition (2) of Theorem~\ref{Thm6.2*a}, that the function $\eta^{1}(\cdot)$ is piecewise constant on $[0,T]$ and admits the representation
\begin{eqnarray}\label{eta}
\eta^{1}(t)=\left\{\begin{array}{ll}
0\quad\mbox{ for a.e. }\;t\in[0,t_1)\;\mbox{ including }t=0,\\
\eta^{1}\;\;\mbox{ for a.e. }\;t\in [t_1,T]\;\mbox{ including }\;t=t_1.
\end{array}\right.
\end{eqnarray}
Since the two robots have the same velocities at the time $t=t_1$ and maintain their velocities until the end of the process, we get $\dot{\ox}^1(t)=\dot{\ox}^2(t)$ for all $t\in[t_1,T]$, which allows us to calculate the value of $\eta^1$ by
\begin{eqnarray}\label{eq-eta2}
\eta^1=\left\{\begin{array}{ll}
\frac{1}{2}\(s_1\ou^1\cos\th_1-s_2\ou^2\cos\th_1\)\;&\mbox{ if }\;s_1\ou^1\ne s_2\ou^2\;\mbox{ and }\;\cos\th_1=\sin\th_1,\\
0&\mbox{ otherwise}.
\end{array}\right.
\end{eqnarray}
Taking into account that the case of $\eta^1=0$ in \eqref{eq-eta2} is trivial, from now on we assume that $\cos\th_1=\sin\th_1$ and $s_1\ou^1\ne s_2\ou^2$.
Remember that after touching the robot pushes the obstacle to the target and they both maintain their constant velocities (speed and direction) until reaching the end of the process at the final time $t=T$. Now using \eqref{nlf}--\eqref{eq-eta2} gives us the trajectory representations
\begin{eqnarray*}
\left\{\begin{array}{ll}
\ox^1(t)=\(\ox^{11}(0),\ox^{12}(0)\)+\(ts_1\ou^1\cos\th,ts_1\ou^1\sin\th_1\),\\
\ox^2(t)=\(\ox^{21}(0),\ox^{22}(0)\)+\(ts_2\ou^2\cos\th_1,ts_2\ou^2\sin\th_1\)\;\mbox{ for }\;t\in[0,t_1),
\end{array}\right.
\end{eqnarray*}\vspace*{-0.2in}
\begin{eqnarray*}
\left\{\begin{array}{ll}
\ox^1(t)=\(\ox^{11}(0),\ox^{12}(0)\)+\big(ts_1\ou^1\cos\th_1-\eta^{1}(t-t_1),ts_1\ou^1\sin\th_1-\eta^{1}(t-t_1)\big),\\
\ox^2(t)=\(\ox^{21}(0),\ox^{22}(0)\)+\big(ts_2\ou^2\cos\th_1+\eta^{1}(t-t_1),ts_2\ou^2\sin\th_1+\eta^{1}(t-t_1)\big)\;\mbox{ for }\;t\in[t_1,T].
\end{array}\right.
\end{eqnarray*}
Employing $\|\ox^2(t_1)-\ox^1(t_1)\|=2R$, we get from the latter formula the equation
\begin{eqnarray}\label{t1}
\begin{array}{ll}
&\[\(s_2\ou^2-s_1\ou^1\)^2\]t_1^2+2\(s_2\ou^2-s_1\ou^1\)\[\(\ox^{21}(0)-\ox^{11}(0)\)\cos\th_1+\(\ox^{22}(0)-\ox^{12}(0)\)\sin\th_1\]t_1\\
&+\(\ox^{21}(0)-\ox^{11}(0)\)^2+\(\ox^{22}(0)-\ox^{12}(0)\)^2-4R^2=0,
\end{array}
\end{eqnarray}
which connects $t_1$ with the given model data and the control $\ou=(\ou^1,\ou^2)$.\vspace*{-0.05in}

To proceed further, for all $t\in[0,T]$ define the functions
\begin{equation*}
d_{ij}(t):=\disp\frac{\ox^i(t)-\ox^j(t)}{\|\ox^i(t)-\ox^j(t)\|},\quad i,j\in\{1,2\}
\end{equation*}
and denote by $\th_{ij}$ the direction of the vector $\ox^i(t_1)-\ox^j(t_1)$. Thus for $t\in[t_1,T]$ we have
\begin{eqnarray*}
\begin{array}{ll}
d_{12}(t)&=\disp\frac{\ox^1(t)-\ox^2(t)}{\|\ox^1(t)-\ox^2(t)\|}=\disp\frac{1}{2R}\(\ox^1(0)-\ox^2(0)+\disp\int_0^{t}\(\dot{\ox}^1(\tau)-\dot{\ox}^2(\tau)\)d\tau\)\\
&=\disp\frac{1}{2R}\(\ox^{11}(0)-\ox^{21}(0)+2t_1\eta^1,\ox^{12}(0)-\ox^{22}(0)+2t_1\eta^1\).
\end{array}
\end{eqnarray*}
On the other hand, it follows from the above that $d_{12}(t_1)=(\cos\th_{12},\sin\th_{12})$, which tells us that
\begin{eqnarray*}
\cos\th_{12}=\frac{\ox^{11}(0)-\ox^{21}(0)+2t_1\eta^1}{2R},\quad\sin\th_{12}=\frac{\ox^{12}(0)-\ox^{22}(0)+2t_1\eta^1}{2R}.
\end{eqnarray*}
This results in determining the value of $y:=t_1\eta^1$ from the quadratic equation
\begin{eqnarray}\label{t.eta}
\begin{array}{ll}
8y^2+4\(x^{11}(0)+x^{12}(0)-x^{21}(0)-x^{22}(0)\)y\\
=4R^2-\((x^{11}(0)-x^{21}(0))^2+(x^{12}(0)-x^{22}(0))^2\).
\end{array}
\end{eqnarray}
Combining \eqref{t.eta} with \eqref{eq-eta2} and \eqref{t1} allows us to precisely compute of optimal solutions when the initial data of the model are specified. The next numerical example illustrates the computation procedure.\vspace{-0.1in}

\begin{example} {\bf (Solving the mobile robot problem without changing direction).}\label{ex1} Specify the model data in the case under consideration by: $n=2,\;x^{01}=\(-30,-30\),\;x^{02}=\(-20,-20\),\;T=6,\;R=6,\;s_1=3,\;s_2=1$ with the compact and control convex set
\begin{equation*}
U:=\big\{u=(u^1,u^2)\in\R^2\big|\;u^1=2u^2,\;-3.37\le u^1\le 3.37\big\}.
\end{equation*}
In this setting we have $t_1>0,\;\th_1=225^{\circ},\;(x^{11}(0)-x^{21}(0))^2+(x^{12}(0)-x^{22}(0))^2=200$. The robot in question has to reach the target by a shortest way, and we assume that the robot tends to maintain its constant direction until either touching the other robot (obstacle), or reaching the end of the process at $t=T$. To proceed with calculations, derive from \eqref{eq-eta2} and \eqref{t.eta} that
\begin{eqnarray*}
\eta^1=\frac{1}{2}\(3\ou^1\(-\frac{\sqrt{2}}{2}\)-\ou^2\(-\frac{\sqrt{2}}{2}\)\)=-\frac{5\sqrt{2}}{4}\ou^2\ne 0\;\mbox{ and }\;
t_1\eta^1=5\pm 3\sqrt{2}.
\end{eqnarray*}
We split our further consideration into the following two cases:\\\vspace*{0.05in}
{\bf Case~1:} $\eta^1=-\frac{5\sqrt{2}}{4}\ou^2$ and $t_1\eta^1=5+3\sqrt{2}$. It gives us $t_1\ou^2=\frac{-12-10\sqrt{2}}{5}$ and the trajectory representations
\begin{eqnarray*}
\left\{\begin{array}{ll}
\ox^1(t)&=\bigg(-25+3\sqrt{2}+\disp\frac{-7\sqrt{2}}{4}\ou^2t,\;-25+3\sqrt{2}+\frac{-7\sqrt{2}}{4}\ou^2t\bigg),\;t\in[t_1,6],\\
\ox^2(t)&=\bigg(-25-3\sqrt{2}+\disp\frac{-7\sqrt{2}}{4}\ou^2t,\;-25-3\sqrt{2}+\frac{-7\sqrt{2}}{4}\ou^2t\bigg),\;t\in[t_1,6].
\end{array}\right.
\end{eqnarray*}
The cost functional is calculated by
\begin{eqnarray*}
J[x,u]=441\left(\ou^2\right)^2+1484.92(\ou^2)+1286
\end{eqnarray*}
and achieves its minimum at $\ou^2\approx-1.68$. Thus $\ou^1\approx-3.37$ and the minimum cost is $J\approx 36$.\\\vspace*{-0.05in}
{\bf Case~2:} $\eta^1=-\frac{5\sqrt{2}}{4}\ou^2$ and $t_1\eta^1=5-3\sqrt{2}$. In this case we get $t_1\ou^2=\frac{12-10\sqrt{2}}{5}$ and
\begin{eqnarray*}
\left\{\begin{array}{ll}
\ox^1(t)&=\bigg(-25-3\sqrt{2}+\disp\frac{-7\sqrt{2}}{4}\ou^2t,\;-25-3\sqrt{2}+\frac{-7\sqrt{2}}{4}\ou^2t\bigg),\;t\in[t_1,6],\\
\ox^2(t)&=\bigg(-25+3\sqrt{2}+\disp\frac{-7\sqrt{2}}{4}\ou^2t,\;-25+3\sqrt{2}+\frac{-7\sqrt{2}}{4}\ou^2t\bigg),\;t\in[t_1,6],
\end{array}\right.
\end{eqnarray*}
with the following expression for the cost functional:
\begin{eqnarray*}
J[x,u]=441\left(\ou^2\right)^2+1484.92(\ou^2)+1286.
\end{eqnarray*}
Thus $J$ achieves its minimum value $J\approx 36$ at $\ou^2\approx-1.68$, and we have $\ou^1\approx-3.37$.\vspace*{-0.05in}

The above calculations show that, in both cases appearing in this setting, the optimal solutions to the robot control problem are calculated as follows:
\begin{eqnarray*}
\left\{\begin{array}{ll}
(\ou^1,\ou^2)&=\(-3.37,-1.68\),\\
\ox^1(t)&=\(-30+7.15t,-30+7.15t\),\;t\in [0,3.11),\\
\ox^1(t)&=\(-20.76+4.16t,-20.76+4.16t\),\;t\in [3.11,6],\\
\ox^2(t)&=\(-20+1.19t,-20+1.19t\),\;t\in [0,3.11),\\
\ox^2(t)&=\(-29.24+4.16t,-29.24+4.16t\),\;t\in [3.11,6].
\end{array}\right.
\end{eqnarray*}\vspace*{-0.05in}

Next we employ the other optimality conditions from Theorem~\ref{Thm6.2*a} to determine adjoint trajectories. Such calculations allow us to reveal more about the optimal model dynamics. It follows from (6) that
\begin{eqnarray*}
\big\la\psi(t),\ou\big\ra=\max_{u\in U}\big\la\psi(t),u\big\ra\mbox{ on }\;[0,6]\;\mbox{ with }\;\psi(t)=\nabla_u g\big(\ox(t),\ou\big)^*q(t),
\end{eqnarray*}
which gives us the equations for the adjoint arc $q(\cdot)$:
\begin{eqnarray*}
s_1\(-\frac{\sqrt{2}}{2}\)q^{11}(t)+s_1\(-\frac{\sqrt{2}}{2}\)q^{12}(t)=\ou^1,\quad s_2\(-\frac{\sqrt{2}}{2}\)q^{21}(t)+s_2\(-\frac{\sqrt{2}}{2}\)q^{22}(t)=\ou^2,
\end{eqnarray*}
and so $q^{11}(t)\approx 0,\;q^{12}(t)\approx 1.59,\;q^{21}(t)\approx 0$, and $q^{22}(t)\approx 2.38$. We deduce from (4) and (7) that $p(t)=p(6)=-\lm\ox(6)-\eta^1x^1_*$ with $\eta^1=-\frac{5\sqrt{2}}{4}\ou^2=2.97$ and $x^1_*=\(1,1,-1,-1\)$. Hence (5) reduces to
\begin{eqnarray*}
\gamma([t,6])=p(t)-q(t)\;\mbox{ on }\;[0,6].
\end{eqnarray*}
Combining it with the above calculations tells us that
\begin{eqnarray*}
\gamma([t,6])=p(6)-q(t)\approx \(-7.17,-7.17,7.25,7.25\)-\(0,1.59,0,2.38\)=\(-7.17,-8.76,7.25,4.87\),
\end{eqnarray*}
for $3.11\le t\le 6$. Thus we confirm that the optimal motion hits the boundary of the state constraint at time $t_1\approx 3.11$ and stays there until the end of the process.
\end{example}\vspace*{-0.15in}

\begin{center}
\begin{figure}[!tbp]
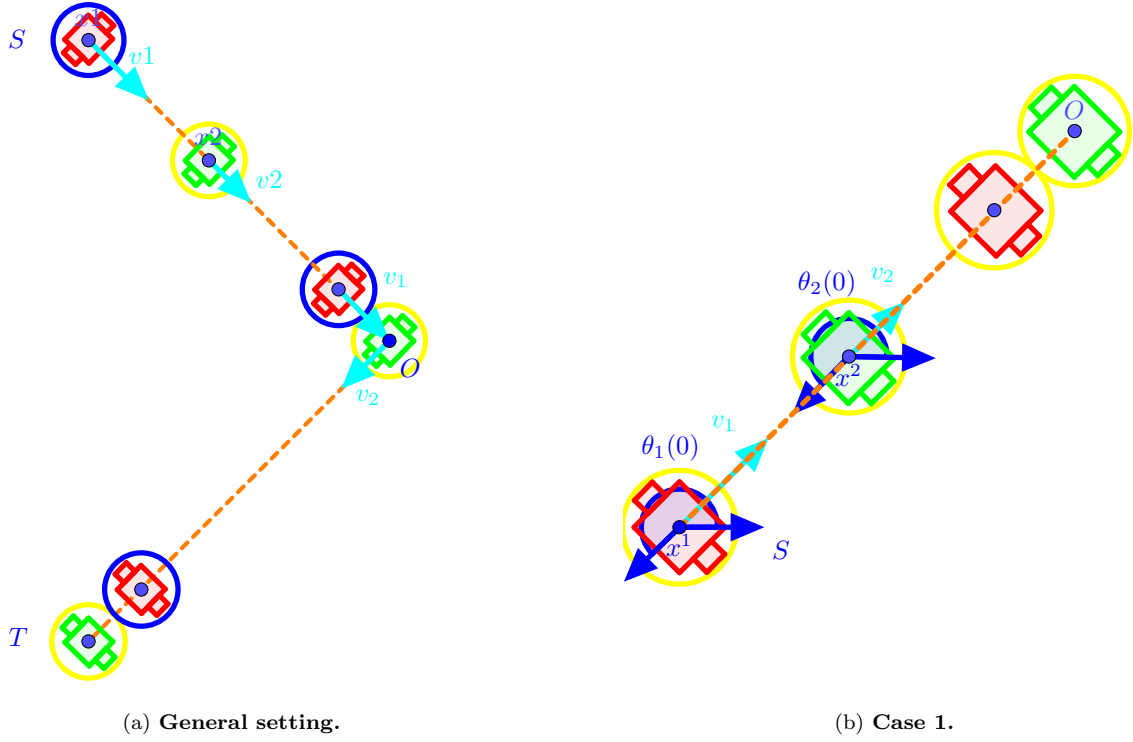

\begin{subfigure}[b]{0.45\textwidth}
\definecolor{qqffff}{rgb}{0.,1.,1.}
\definecolor{qqffqq}{rgb}{0.,1.,0.}
\definecolor{ffqqqq}{rgb}{1.,0.,0.}
\definecolor{ffxfqq}{rgb}{1.,0.4980392156862745,0.}
\definecolor{ududff}{rgb}{0.30196078431372547,0.30196078431372547,1.}
\definecolor{ffffqq}{rgb}{1.,1.,0.}
\definecolor{qqqqff}{rgb}{0.,0.,1.}

\vspace*{-0.1in}
\caption{\bf General setting.}
\label{picnew}
\end{subfigure}
\hfill
\begin{subfigure}[b]{0.45\textwidth}
\definecolor{ffxfqq}{rgb}{1.,0.4980392156862745,0.}
\definecolor{qqffff}{rgb}{0.,1.,1.}
\definecolor{qqffqq}{rgb}{0.,1.,0.}
\definecolor{ffqqqq}{rgb}{1.,0.,0.}
\definecolor{xdxdff}{rgb}{0.49019607843137253,0.49019607843137253,1.}
\definecolor{qqqqff}{rgb}{0.,0.,1.}
\definecolor{ffffqq}{rgb}{1.,1.,0.}
\definecolor{ududff}{rgb}{0.30196078431372547,0.30196078431372547,1.}
%
\vspace*{-0.1in}
\caption{\bf Case~1.}
\label{pic2a*}
\end{subfigure}
\caption{\bf Mobile robot model with changing directions in contact.}
\end{figure}
\end{center}

$\textbf{B: Mobile robot model with changing direction in contact.}$ Now we examine other situations in robot behavior before and after contacting the obstacle that are different from the previous consideration in setting A. Let $t_1$ be the contacting time as in \eqref{t-contact}. Consider the case where robot~1 in question moves faster than robot~2 (obstacle) and touches the second robot at $t_1$, while after the contact both robots together change their directions to go to the target with the same speed; see Fig.~\ref{picnew}. We have\vspace*{-0.45in}
\begin{eqnarray*}
\th_1(t)=\th_2(t)=\left\{\begin{array}{ll}
\th_1(0)\;\mbox{ if }\;t\in [0,t_1),\\
\th_1(t_1)\;\mbox{ if }\;t\in [t_1,T],
\end{array}\right.
\end{eqnarray*}
where $\th_j(t)$, $j=1,2$, are angles of the corresponding robot directions. Before the time $t_1$ both robots move in the same direction with different speeds, but at the contact time $t_1$ they change their directions and go together to the target with the same speed. Thus the velocities of the two robots are given by
\begin{eqnarray*}
\left\{\begin{array}{ll}
\dot{\ox}^1(t)=\(s_1\ou^1\cos\th_1(0),\;s_1\ou^1\sin\th_1(0)\),\\
\dot{\ox}^2(t)=\(s_2\ou^2\cos\th_1(0),\;s_2\ou^2\sin\th_1(0)\)\;\mbox{ for }\;t\in[0,t_1),
\end{array}\right.
\end{eqnarray*}\vspace*{-0.2in}
\begin{eqnarray*}
\left\{\begin{array}{ll}
\dot{\ox}^1(t)=\(s_1\ou^1\cos\th_1(t_1)-\eta^{1}(t),\;s_1\ou^1\sin\th_1(t_1)-\eta^{1}(t)\),\\
\dot{\ox}^2(t)=\(s_2\ou^2\cos\th_1(t_1)+\eta^{1}(t),\;s_2\ou^2\sin\th_1(t_1)+\eta^{1}(t)\)\;\mbox{ for }\;t\in[t_1,T],
\end{array}\right.
\end{eqnarray*}
where the piecewise constant function $\eta^{1}(\cdot)$ on $[0,T]$ is taken from \eqref{eta}. Furthermore, starting with the contact time $t=t_1$ the robots tend to maintain the same velocities until the end of the process. This allows us to calculate the value of $\eta^1$ in \eqref{eta} by
\begin{eqnarray}\label{eq-eta2-game}
\eta^1=\left\{\begin{array}{ll}
\frac{1}{2}\(s_1\ou^1\cos\th_1(t_1)-s_2\ou^2\cos\th_1(t_1)\)\;&\mbox{ if }\;s_1\ou^1\ne s_2\ou^2\;\mbox{ and }\;\cos\th_1(t_1)=\sin\th_1(t_1),\\
0&\mbox{ otherwise}.
\end{array}\right.
\end{eqnarray}
Excluding the trivial case $\eta^1=0$, suppose that $\cos\th_1(t_1)=\sin\th_1(t_1)$ and $s_1\ou^1\ne s_2\ou^2$. Similarly to our previous consideration in A, we arrive at the same quadratic equation \eqref{t.eta} for the value $y:=t_1\eta_1$, but in the new setting. The corresponding trajectory representations are given now by
\begin{eqnarray*}
\left\{\begin{array}{ll}
\ox^1(t)=\(\ox^{11}(0)+s_1\ou^1\cos\th_1(0)t,\;\ox^{12}(0)+s_1\ou^1\sin\th_1(0)t\),\\
\ox^2(t)=\(\ox^{21}(0)+s_2\ou^2\cos\th_1(0)t,\;\ox^{22}(0)+s_2\ou^2\sin\th_1(0)t\)\;\mbox{ for }\;t\in[0,t_1),
\end{array}\right.
\end{eqnarray*}\vspace*{-0.2in}
\begin{eqnarray*}
\left\{\begin{array}{ll}
\ox^1(t)&=\big(\ox^{11}(0)+s_1\ou^1\cos\th_1(0)t_1+s_1\ou^1\cos\th_1(t_1)(t-t_1)-\eta^1(t-t_1),\\
&\quad\ox^{12}(0)+s_1\ou^1\sin\th_1(0)t_1+s_1\ou^1\sin\th_1(t_1)(t-t_1)-\eta^1(t-t_1)\big),\\
\ox^2(t)&=\big(\ox^{21}(0)+s_2\ou^2\cos\th_1(0)t_1+s_2\ou^2\cos\th_1(t_1)(t-t_1)+\eta^1(t-t_1),\\
&\quad\ox^{22}(0)+s_2\ou^2\sin\th_1(0)t_1+s_2\ou^2\sin\th_1(t_1)(t-t_1)+\eta^1(t-t_1)\big)\;\mbox{ for }\;t\in[t_1,T].
\end{array}\right.
\end{eqnarray*}
Combining the noncollision conditions with \eqref{eq-eta2-game} allows us to conclude that
\begin{eqnarray}\label{restrictions}
\left\{\begin{array}{ll}
&x^{21}(t)>x^{11}(t)\mbox{ and }x^{22}(t)>x^{12}(t)\;\mbox{ for all }\;t\in[0,T],\\
&\mbox{with }\;s_1\ou^1\ne s_2\ou^2\mbox{ and }\cos\th_1(t_1)=\sin\th_1(t_1),
\end{array}\right.
\end{eqnarray}
which lead us to the following two cases.\\[1ex]
$\bullet$ {\bf Case~1 (robots are in the third quadrant):} In this case illustrated by Fig.~\ref{pic2a*} we have
\begin{eqnarray*}
\th_1(t)=\th_2(t)=:\th\;\mbox{ for all }\;t\in [0,T].
\end{eqnarray*}
Prior to the final time $t_1=T$ both robots move in the same direction with different speeds, while at $t_1=T$ they are in contact and reach the target.
Thus we get from \eqref{restrictions} that
\begin{eqnarray}\label{restrictions*}
\left\{\begin{array}{ll}
&x^{21}(t)>x^{11}(t)\mbox{ and }x^{22}(t)>x^{12}(t)\;\mbox{ for all }\;t\in[0,T]\\
&\mbox{with }\;s_1\ou^1\ne s_2\ou^2\;\mbox{ and }\;\cos\th=\sin\th,
\end{array}\right.
\end{eqnarray}
and that all the other formulas above hold with the corresponding specifications. Note that in this case both robots reach the target at the final time $t=T$, and the minimum cost is $J=0$. It obviously shows that the distance from the robot to the target is the shortest one.\vspace*{-0.1in}
\begin{center}
\begin{figure}[!tbp]
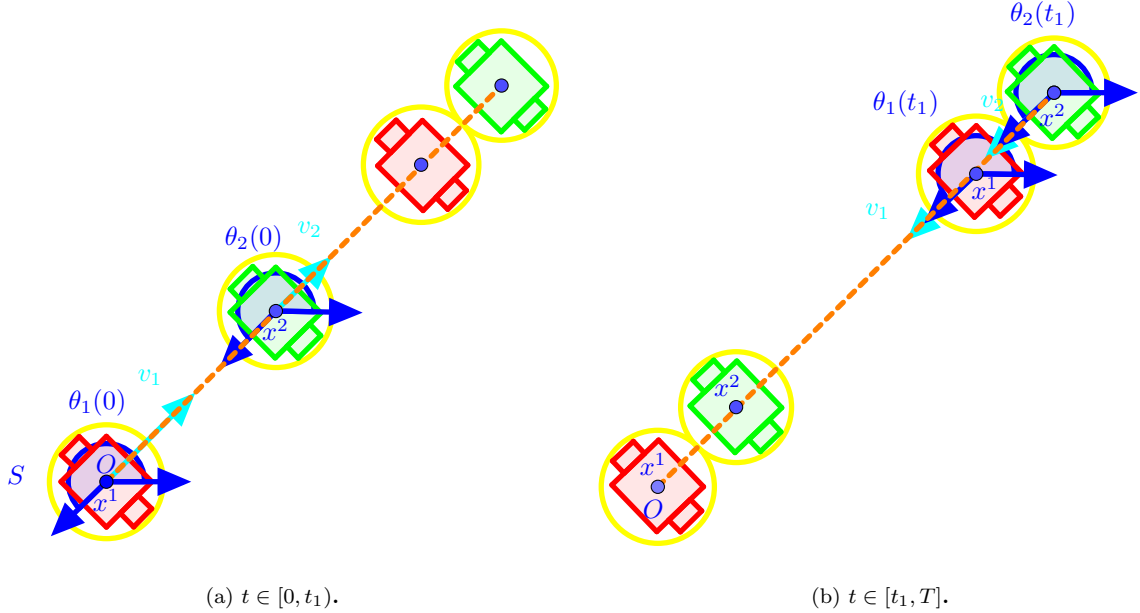

\begin{subfigure}[b]{0.5\textwidth}
\definecolor{ffxfqq}{rgb}{1.,0.4980392156862745,0.}
\definecolor{qqffff}{rgb}{0.,1.,1.}
\definecolor{qqffqq}{rgb}{0.,1.,0.}
\definecolor{ffqqqq}{rgb}{1.,0.,0.}
\definecolor{xdxdff}{rgb}{0.49019607843137253,0.49019607843137253,1.}
\definecolor{qqqqff}{rgb}{0.,0.,1.}
\definecolor{ffffqq}{rgb}{1.,1.,0.}
\definecolor{ududff}{rgb}{0.30196078431372547,0.30196078431372547,1.}

\vspace*{-0.2in}
\caption{\bf $t\in [0,t_1)$.}
\label{pic2a}
\end{subfigure}
\hfill
\begin{subfigure}[b]{0.5\textwidth}
\definecolor{ffxfqq}{rgb}{1.,0.4980392156862745,0.}
\definecolor{qqffff}{rgb}{0.,1.,1.}
\definecolor{qqffqq}{rgb}{0.,1.,0.}
\definecolor{ffqqqq}{rgb}{1.,0.,0.}
\definecolor{xdxdff}{rgb}{0.49019607843137253,0.49019607843137253,1.}
\definecolor{qqqqff}{rgb}{0.,0.,1.}
\definecolor{ffffqq}{rgb}{1.,1.,0.}
\definecolor{ududff}{rgb}{0.30196078431372547,0.30196078431372547,1.}
%
\vspace*{-0.2in}
\caption{\bf $t\in[t_1,T]$.}
\label{pic2b}
\end{subfigure}
\caption{\bf Mobile robot model with changing directions in Case~2.}
\end{figure}
\end{center}\vspace*{-0.3in}
$\bullet$ {\bf Case~2 (robots are in the first quadrant):} In this case robot~1 in question moves faster than the robot~2 and touches the latter at the contact time $t_1\ne T$. Then robot~1 pulls robot~2 to go back to the starting point with the same speed, where the starting point is taken as the target at the origin. Then we also have $\th_1(t)=\th_2(t)=\th\;\mbox{ for all }\;t\in[0,T]$; see Fig.~\ref{pic2a} and Fig.~\ref{pic2b}. Prior to the contact time $t_1$ both robots move in the same direction with different speeds, while at the contact time $t_1$ they change their directions simultaneously and move together to the starting point with the same speed. Thus we can proceed similarly to Case~1 under the conditions in \eqref{restrictions*}.\vspace*{-0.25in}

\section{Controlled Model of Pedestrian Traffic Flows}\vspace*{-0.1in}

In this section we formulate a continuous-time, deterministic, and optimal control version of the pedestrian traffic flow model through a doorway for which a stochastic, discrete-time, and simulation (uncontrolled) counterpart was originated in \cite{GG}. Here we formalize the dynamics via a perturbed sweeping process with constrained controls in perturbations that should be determined to ensure the desired performance. We also discuss differences and similarities with the crowd motion model of the pedestrian traffic as well as with the mobile robot model formulated and studied in Section~3.\vspace*{-0.05in}

In the model under consideration we have $n$ pedestrians $x^i\in\R,\;i=1,\ldots,n$ as $n\ge 2$ that are identified with rigid disks of the same radius $R$ going through a doorway as depicted in Fig.~\ref{pic3}.\vspace*{-0.05in}
\begin{center}
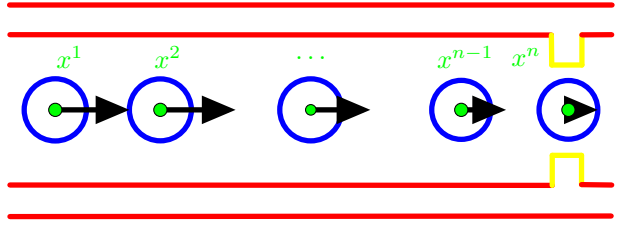
\begin{figure}[htp]
\begin{center}
\definecolor{ffffqq}{rgb}{1.,1.,0.}
\definecolor{qqqqff}{rgb}{0.,0.,1.}
\definecolor{qqffqq}{rgb}{0.,1.,0.}
\definecolor{ffqqqq}{rgb}{1.,0.,0.}
\begin{tikzpicture}[line cap=round,line join=round,>=triangle 45,x=1.0cm,y=1.0cm]
\clip(-4.766730466479051,-2.1545861011508163) rectangle (4.875864531668342,2.1948397272137834);
\draw [line width=2.pt,color=ffqqqq] (-4.004577708577143,1.3950042436343764)-- (4.004432685958349,1.3950042436343764);
\draw [line width=2.pt,color=ffqqqq] (-4.,-1.42)-- (3.9926373834627156,-1.412277750326518);
\draw [line width=2.pt,color=qqqqff] (-3.3969044328607474,0.) circle (0.4135733067081072cm);
\draw [line width=2.pt,color=qqqqff] (-2.,0.) circle (0.4118734884757825cm);
\draw [line width=2.pt,color=qqqqff] (0.,0.) circle (0.41257379322447274cm);
\draw [line width=2.pt,color=qqqqff] (2.,0.) circle (0.39239240244126683cm);
\draw [line width=2.pt,color=qqqqff] (3.424615575509039,0.) circle (0.39031032939731797cm);
\draw [line width=2.pt,color=ffqqqq] (-4.,1.)-- (3.198315720112725,0.9952165148262837);
\draw [line width=2.pt,color=ffffqq] (3.198315720112725,0.9952165148262837)-- (3.203048520883285,0.5967973784567968);
\draw [line width=2.pt,color=ffffqq] (3.203048520883285,0.5967973784567968)-- (3.604162900550126,0.5967973784567968);
\draw [line width=2.pt,color=ffffqq] (3.604162900550126,0.5967973784567968)-- (3.604162900550126,0.9979117581236371);
\draw [line width=2.pt,color=ffqqqq] (3.604162900550126,0.9979117581236371)-- (4.,1.);
\draw [line width=2.pt,color=ffqqqq] (-4.,-1.)-- (3.2095817077789257,-1.004840635613907);
\draw [line width=2.pt,color=ffffqq] (3.2095817077789257,-1.004840635613907)-- (3.2095817077789257,-0.6050956690227378);
\draw [line width=2.pt,color=ffffqq] (3.2095817077789257,-0.6050956690227378)-- (3.599076803431861,-0.6050956690227378);
\draw [line width=2.pt,color=ffffqq] (3.599076803431861,-0.6050956690227378)-- (3.599076803431861,-0.9945907646756719);
\draw [line width=2.pt,color=ffqqqq] (3.599076803431861,-0.9945907646756719)-- (4.,-1.);
\draw [->,line width=2.pt] (3.424615575509039,0.) -- (3.8149259049063557,0.);
\draw [->,line width=2.pt] (2.,0.) -- (2.5958161381638964,0.);
\draw [->,line width=2.pt] (0.,0.) -- (0.7916408691171811,0.);
\draw [->,line width=2.pt] (-2.,0.) -- (-1.,0.);
\draw [->,line width=2.pt] (-3.3969044328607474,0.) -- (-2.4118734884757824,0.);
\begin{normalsize}
\draw [fill=qqffqq] (-3.3969044328607474,0.) circle (2.5pt);
\draw[color=qqffqq] (-3.2,0.7) node {$x^1$};
\draw [fill=qqffqq] (-2.,0.) circle (2.5pt);
\draw[color=qqffqq] (-1.9,0.7) node {$x^2$};
\draw [fill=qqffqq] (0.,0.) circle (2.0pt);
\draw[color=qqffqq] (0,0.7) node {$\ldots$};
\draw [fill=qqffqq] (2.,0.) circle (2.5pt);
\draw[color=qqffqq] (2.0548925907209665,0.7) node {$x^{n-1}$};
\draw [fill=qqffqq] (3.424615575509039,0.) circle (2.5pt);
\draw[color=qqffqq] (2.865280893756758,0.7) node {$x^n$};
\end{normalsize}
\end{tikzpicture}
\end{center}
\caption{Unidirectional flows of pedestrians through doorway.}
\label{pic3}
\end{figure}
\end{center}\vspace*{-0.1in}

Define the set of {\em admissible configurations} by imposing the {\em nonoverlapping conditions} in order to avoid overlapping between two pedestrians:
\begin{eqnarray}\label{e:129}
Q_0:=\big\{x=\big(x^1,\ldots,x^n\big)\in\R^n\big|\;x^{i+1}-x^i\ge 2R\;\mbox{ whenever }\;i,j\in\{1,\ldots,n\}\big\}.
\end{eqnarray}
Denoting by $S(x)$ the {\em spontaneous velocity} of the pedestrians at $x\in Q_0$, we represent it as
\begin{eqnarray*}
S(x):=\big(S_0(x^1),\ldots,S_0(x^n)\big)\;\mbox{ with }\;S_0(x)=s_0\nabla D(x),\quad x\in Q_0,
\end{eqnarray*}
where $Q_0$ is taken from \eqref{e:129}, $D(x)$ denotes the {\em distance} from the position $x=(x^1,\ldots,x^n)\in Q_0$ to the doorway, and the scalar $s_0\ge 0$ indicates the speed. Since $x\ne 0$ and hence $\|\nabla D(x)\|=1$, we get $s_0=\|S_0(x)\|$. Each pedestrian tends to maintain his/her desired spontaneous velocity until reaching the doorway in the absence of other pedestrians that is reflected in the model by
\begin{eqnarray}\label{g-control}
g(x)=(s_1,\ldots,s_n)\in\R^n\;\textrm{ for all }\;x=(x_1,\ldots,x_n)\;\in Q_0,
\end{eqnarray}
where $s_i$ denotes the speed of the pedestrian $i\in\{1,\ldots,n\}$. If the distance between pedestrian $i$ and pedestrian $i+1$ is $x^{i+1}(t)-x^i(t)=2R$, then both pedestrians tend to adjust their velocities in order to keep the distance to be at least $2R$. In this setting we use some force in order to control the actual velocity of all the pedestrians in the presence of the nonoverlapping conditions \eqref{e:129}. This is modeled by inserting controls $u(\cdot)=(u^1(\cdot),\ldots,u^n(\cdot))$ into the perturbation term as follows:
\begin{eqnarray}\label{cr-pert}
g(x(t),u(t)):=\big(s_1u^1(t),\ldots,s_nu^n(t)\big),\quad t\in[0,T],
\end{eqnarray}
where measurable control functions $u=(u^1,\ldots,u^n)\colon[0,T]\to\R^n$ satisfy the constraint
\begin{eqnarray}\label{cr-u}
u(t)\in U\;\textrm{ a.e. on }\;[0,T]
\end{eqnarray}
defined via a convex and compact set $U\subset\R^n$, which is specified below in particular situations.

Observing that the pedestrians cannot move with their spontaneous velocities due to the nonoverlapping constraints in \eqref{e:129}, we consider the set of {\em feasible velocities}
\begin{eqnarray*}
V_x:=\big\{v=\big(v^1,\ldots,v^n\big)\in\R^n\left|\right.x^{i+1}-x^i=2R\Longrightarrow v^{i+1}\ge v^i\;\mbox{ for all }\;i=1,\ldots,n-1\big\}
\end{eqnarray*}
and then describe the {\em actual velocity field} is the feasible field in terms of the (unique) Euclidean projection of the spontaneous velocity $S\big(x\big)$ onto the convex set $V_x$ by
\begin{eqnarray}\label{desc}
\dot{x}(t)=\Pi\big(S(x);V_x\big)\;\mbox{ for a.e. }\;t\in[0,T],\quad x(0)=x_0\in Q_0,
\end{eqnarray}
where $x_0$ indicates the starting position of the pedestrians. Based on the projection description \eqref{desc} and definition \eqref{1.4} of the normal cone of convex analysis, we deduce from \eqref{desc} that
\begin{eqnarray*}
S\big(x\big)\in N(x;Q_0)+\dot{x}(t)\;\mbox{ for a.e. }\;t\in[0,T],\quad x(0)=x_0,
\end{eqnarray*}
which gives us the differential inclusion of the perturbed sweeping process
\begin{eqnarray}\label{sw-desc}
\dot{x}(t)\in-N(x;Q_0)+S(x)\;\mbox{ for a.e. }\;t\in[0,T],\quad x(0)=x_0.
\end{eqnarray}
Define further the convex set $C\subset\R^n$ by
\begin{eqnarray}\label{e:131*}
C:=\big\{x\in\R^n\big|\;\la x^j_*,x\ra\le c_j,\;j=1,\ldots,n-1\big\}\;\mbox{ with }\;x^j_*:=e_j-e_{j+1},\;c_j=-2R
\end{eqnarray}
for $j=1,\ldots,n-1$, where $(e_1,\ldots,e_n)$ are the orths in $\R^n$. Remembering the control velocity description \eqref{g-control} allows us to describe the pedestrian model dynamics as the controlled sweeping process
\begin{eqnarray}\label{e:135*}
\left\{\begin{array}{lcl}
-\dot{x}(t)\in N\big(x(t);C\big)-g\(x(t),u(t)\)\;\mbox{ for a.e. }\;t\in[0,T],\\
u(t)\in U\;\mbox{ for a.e. }\;t\in[0,T],\quad x(0)=x_0\in C,
\end{array}\right.
\end{eqnarray}
with $C$ and $U$ taken from \eqref{e:131*} and \eqref{cr-u}, respectively. Note the differential inclusion in \eqref{e:135*} intrincically contains the pointwise state constraint
\begin{eqnarray}\label{cr-state}
x(t)\in C\;\mbox{ for all }\;t\in[0,T],
\end{eqnarray}
which is equivalent to the nonoverlapping conditions from \eqref{e:129} due to the structure of $C$ in \eqref{e:131*}.\vspace*{-0.05in}

Furthermore, it surely makes sense to introduce an appropriate {\em cost functional} to optimize the performance of the model over the constrained dynamics in \eqref{e:135*} and to formulate an optimal control problem in the form of $(P)$ from Section~2. A very natural candidate for the cost functional, which reflects the essence and goal of the model, is the following one:
\begin{eqnarray}\label{cr-cost}
\mbox{minimize }\;J[x,u]:=\frac{1}{2}\|x(T)\|^2
\end{eqnarray}
meaning the minimization of the distance from all the pedestrians from \eqref{e:129} to the doorway at the origin. The obtained description of the controlled pedestrian traffic model allows us applying the necessary optimality conditions for problem $(P)$ presented in Section~2 to find optimal solutions in this model that exist due to \cite{et}. Prior to such an application, let us compare the model under consideration with those for the controlled crown model from \cite{cm2} and for the mobile robot model studied in Section~3.\vspace*{-0.15in}

\begin{remark}{\bf(comparison with the crowd motion model).} There are certain similarities between the controlled pedestrian traffic flow model through a doorway considered here and the optimization model for controlled crowd motions in a corridor studied in \cite{cm2} via alternative necessary optimality conditions for absolutely continuous controls of a perturbed sweeping process. However, a crucial difference of the present model from the one considered in \cite{cm2} is that now we are able, based on the new results of \cite{cmn}, to deal with real-life pointwise constraints on control functions, which are unavoidable in practice while being highly theoretically challenging. Incorporating such constraints allows us to exclude the energy term from the cost functional and concentrate on minimizing the distance of participants from the target, which adequately reflects the very essence of the model.\vspace*{-0.05in}

Mathematically we can treat the pointwise (hard) control constraints by the powerful {\em maximum principle} established in \cite{cmn} for the controlled perturbed sweeping process under consideration; see more details below. This was not the case in the unconstrained setting of \cite{cm2}.
\end{remark} \vspace*{-0.25in}

\begin{remark}{\bf(comparison with the mobile robot model).} Although the essence and practical sense of the controlled robotics model studied in Section~3 and the controlled pedestrian flow model considered in this section are completely different, there are some similarities in their mathematical descriptions as perturbed sweeping processes. This allows us to apply the same necessary optimality conditions from \cite{cmn} to determining optimal solutions in both models. Of course, the main mathematical difference between the mathematical descriptions of these two models is the space dimension. On the other hand, the available results for planar crowd motion models developed in \cite{cm3,cm4} are not applicable to either of the models considered in Sections~3 and 4 due to the unconstrained nature of the previously obtained developments.
\end{remark}\vspace*{-0.1in}

Denoting now by $(SF)$ the optimal control problem for the pedestrian traffic flow model formulated by \eqref{cr-state}--\eqref{cr-cost} with the data from \eqref{cr-pert}--\eqref{e:131*}, we apply to it the necessary optimality conditions of Theorem~\ref{Thm6.2*} obtained for problem $(P)$ of this category. The next theorem specifies the obtained results in the case of problem $(SF)$ under consideration.\vspace*{-0.1in}

\begin{theorem}\label{Thm6.2*c}{\bf (necessary optimality conditions for the sweeping control pedestrian traffic flow model).} Let $(\ox(\cdot),\ou(\cdot))$ be a $W^{1,2}\times L^2$-local minimizer of problem $(SF)$, where the control set $U$ is compact and convex. Then there exist a multiplier $\lm\ge 0$, a measure $\gg=(\gg^1,\ldots,\gg^n)\in C^*([0,T];\R^n)$ as well as adjoint arcs $p(\cdot)\in W^{1,2}([0,T];\R^n)$ and $q(\cdot)\in BV([0,T];\R^n)$ satisfying to the following conditions:\vspace*{-0.1in}
\begin{itemize}
\item[\bf(1)] $-\dot{\ox}(t)=\disp\sum^{n-1}_{j=1}\eta^j(t)x_*^j-\big(s_1\ou^1(t),\ldots,s_n\ou^n(t)\big)$ for a.e.\ $t\in[0,T]$,\\where $\eta^j(\cdot)
\in L^2([0,T];\R_+)$ are uniquely defined by this representation and well defined at $t=T$;
\item[\bf(2)] $\ox^{j+1}(t)-\ox^j(t)>2R\Longrightarrow\eta^j(t)=0$ for all $j=1,\ldots,n-1$ and a.e. $t\in[0,T]$ including $t=T$;
\item[\bf(3)] $\eta^j(t)>0\Longrightarrow\la x^j_*,q(t)\ra=c_j$ for all $j=1,\ldots,n-1$ and a.e. $t\in[0,T]$ including $t=T$;
\item[\bf(4)] $p(t)=p(T)$ for all $t\in[0,T]$;
\item[\bf(5)] $q(t)=p(T)-\gamma([t,T])$ for all $t\in[0,T]$ except at most a countable subset;
\item[\bf(6)] $\big\la\psi(t),\ou(t)\big\ra=\max_{u\in U}\big\la\psi(t),u\big\ra\;\mbox{ for a.e. }\;t\in[0,T]$, where $\psi(t):=\begin{pmatrix}
s_1&0&\ldots&0\\
0&s_2&\ldots&0\\
\ldots&\ldots&\ldots&\ldots\\
0&0&\ldots&s_n
\end{pmatrix}q(t)$;
\item[\bf(7)] $-p(T)=\lm\ox(T)+\sum_{j\in I(\ox(T))}\eta^j(T)x^j_*$ via the set of active constraint indices $I(\ox(T))$ at $\ox(T)$;
\item[\bf(8)] $\sum_{j\in I(\ox(T))}\eta^j(T)x^j_*\in N\big(\ox(T);C)$;
\item[\bf(9)] $\lm+\|q(0)\|+\|p(T)\|>0$.
\end{itemize}
\end{theorem}\vspace*{-0.1in}
{\bf Proof.} It is direct consequence of Theorem~\ref{Thm6.2*} with the data of $(P)$ specified for $(SF)$ by particular taking into account the form of the controlled perturbation mapping $g$ in \eqref{cr-pert}. $\h$

Let us discuss some immediate conclusions for the pedestrian traffic flow model that can be derived from the obtained theorem.\vspace*{-0.05in}

$\bullet$ At the contacting time $t_1\in[0,T]$ when $\ox^{i+1}(t_1)-\ox^i(t_1)=2R,\;i=1,\ldots,n-1$, pedestrians $i$ and $i+1$ adjust their speeds in order to keep the distance between them to be at least $2R$. It is natural to suppose that after the time $t=t_1$ both pedestrians $i$ and $i+1$ tend to maintain their new constant velocities until either reaching someone ahead or stopping at $t=T$. Hence the velocities of all the pedestrians are piecewise constant on $[0,T]$ in this setting.\vspace*{-0.05in}

$\bullet$  The controlled system of the differential equations in (1) can be written as
\begin{eqnarray}\label{e:135a}
\left\{\begin{array}{ll}
-\dot{\ox}^1(t)=\eta^1(t)-s_1\ou^1(t),\\
-\dot{\ox}^i(t)=\eta^i(t)-\eta^{i-1}(t)-s_i\ou^i(t),\quad i=2,\ldots,n-1,\\
-\dot{\ox}^n(t)=-\eta^{n-1}(t)-s_n\ou^n(t)\;\mbox{ for a.e. }\;t\in[0,T].
\end{array}\right.
\end{eqnarray}
If pedestrian~1 does not touch pedestrian~2 in the sense that $\ox^{2}(t)-\ox^1(t)>2R$ for all $t\in[0,T]$, then it follows from \eqref{e:135a} and (2) that the actual velocity and the spontaneous velocity of pedestrian~1 agree for a.e.\ $t\in[0,T]$, which means that $\dot{\ox}^1(t)=s_1\ou^1(t)$ a.e.\ on $[0,T]$. If $\ox^{n}(t)-\ox^{n-1}(t)>2R$ for all $t\in[0,T]$, we get this conclusion for pedestrian~$n$. The same holds for pedestrians $i=2,\ldots,n-1$ provided that $\ox^{i+1}(t)-\ox^i(t)>2R$ and $\ox^{i}(t)-\ox^{i-1}(t)>2R$ whenever $t\in[0,T]$.\vspace*{-0.05in}

To proceed further, suppose that $\lm>0$ (say $\lm=1$); otherwise, it is not enough information to efficiently apply Theorem~\ref{Thm6.2*a}. Moreover, assuming for simplicity of calculations in the examples below that the control actions $\ou^i(\cdot)$ are constant $\ou^i$ on $[0,T]$ for all $i=1,\ldots,n$ and then employing the Newton-Leibniz formula in \eqref{e:135a} gives us the trajectories
\begin{eqnarray}\label{e:137*}
\left\{\begin{array}{ll}
\ox^1(t)=x^{01}-\disp\int^t_0\eta^1(\tau)d\tau+ts_1\ou^1,\\
\ox^i(t)=x^{0i}+\disp\int^t_0\big[\eta^{i-1}(\tau)-\eta^i(\tau)\big]d\tau+ts_i\ou^i\;\mbox{ as }\;i=2,\ldots,n-1,\\
\ox^n(t)=x^{0n}+\disp\int^t_0\eta^{n-1}(\tau)d\tau+ts_n\ou^n
\end{array}\right.
\end{eqnarray}
for all $t\in[0,T]$, where $(x^{01},\ldots,x^{0n})$ are the components of the starting point $x_0\in C$ in \eqref{e:135*}.

Next we fix $i\in\{1,\ldots,n-1\}$ and let $t_i$ be the first time when $\ox^{i+1}(t_i)-\ox^i(t_i)=2R$. Observe that the the vector function $\eta(\cdot)$ in the conditions above is piecewise constant on $[0,T]$ and rewrite \eqref{e:137*} by
\begin{eqnarray*}
\ox^i(t)=x^{0i}+\disp{\int^t_0}\big[\eta^{i-1}(\tau)-\eta^i(\tau)\big]d\tau+ts_i\ou^i\;\mbox{ for }\;i=1,\ldots,n
\end{eqnarray*}
with $\eta^0=\eta^n=0$. For each $i$ define the positive numbers $\Theta^i$ and $\Theta_i$ by
\begin{eqnarray*}
\Theta^i:=\min\big\{t_j\big|\;t_j>t_i,\;j=1,\ldots,n-1\big\},\quad\Theta_i:=\max\big\{t_j\big|\;t_j<t_i,\;j=1,\ldots,n-1\big\}.
\end{eqnarray*}
Then we have the following trajectory representations:
\begin{eqnarray*}
\left\{\begin{array}{ll}
\ox^i(t)=x^{0i}+\disp{\int^t_0}\eta^{i-1}(\tau)d\tau+ts_i\ou^i,\quad t\in[0,t_i),\\
\ox^i(t)=x^{0i}+\disp\int^{t_i}_0\eta^{i-1}(\tau)d\tau+(t-t_i)\big[\eta^{i-1}(t_i)-\eta^i(t_i)\big]+ts_i\ou^i,\quad t\in[t_i,\Theta^i],
\end{array}\right.
\end{eqnarray*}\vspace*{-0.15in}
\begin{eqnarray*}
\left\{\begin{array}{ll}
\ox^{i+1}(t)=x^{0(i+1)}-\disp{\int^t_0}\eta^{i+1}(\tau)d\tau+ts_{i+1}\ou^{i+1},\quad t\in[0,t_i),\\
x^{i+1}(t)=x^{0(i+1)}-\disp\int^{t_i}_0\eta^{i+1}(\tau)d\tau+(t-t_i)\big[\eta^i(t_i)-\eta^{i+1}(t_i)\big]+ts_{i+1}\ou^{i+1},\quad t\in[t_i,\Theta^i].
\end{array}\right.
\end{eqnarray*}
Suppose without loss of generality that the functions $\dot\ox^i(\cdot)$ are well defined at $t_i$ while $\eta^i(\cdot)$ are well defined at $t_i$ and $\Theta_i$. At the contact time $t=t_i$ we get $\ox^{i+1}(t_i)-\ox^i(t_i)=2R$ and
\begin{eqnarray*}
&\;&\ox^{i+1}(t_i)-\ox^i(t_i)\\
&=&x^{0(i+1)}-x^{0i}-\int^{t_i}_0\big[\eta^{i+1}(\tau)+\eta^{i-1}(\tau)\big]d\tau+t_i\big(s_{i+1}\ou^{i+1}-s_i\ou^i\big)\\
&=&x^{0(i+1)}-x^{0i}-\int^{\Theta_i}_0\big[\eta^{i+1}(\tau)+\eta^{i-1}(\tau)\big]d\tau-(t_i-\Theta_i)\big[\eta^{i+1}(\Theta_i)+\eta^{i-1}(\Theta_i)\big]
+t_i(s_{i+1}\ou^{i+1}-s_i\ou^i).
\end{eqnarray*}
Then we arrive at the following conclusions:\vspace*{-0.05in}

$\bullet$ If $x^{0(i+1)}-x^{0i}=2R$, it is easy to see that $t_i=0$.\vspace*{-0.05in}

$\bullet$ If $x^{0(i+1)}-x^{0i}>2R$, it follows
that
\begin{eqnarray}\label{e:145}
t_i=\dfrac{x^{0(i+1)}-x^{0i}-2R+\Theta_i\big[\eta^{i+1}(\Theta_i)+\eta^{i-1}(\Theta_i)\big]-\disp{\int^{\Theta_i}_0}\big[\eta^{i+1}(\tau)+\eta^{i-1}
(\tau)\big]d\tau}{\eta^{i+1}(\Theta_i)+\eta^{i-1}(\Theta_i)-s_{i+1}\ou^{i+1}+s_i\ou^i}.
\end{eqnarray}
Since after the contact at $t_i$ the pedestrians go to the target with the same velocity, we get
\begin{eqnarray}\label{e:147}
\dot{\ox}^{i+1}(t_i)=\dot{\ox}^{i}(t_i)\Longleftrightarrow 2\eta^i(t_i)=\eta^{i+1}(t_i)+\eta^{i-1}(t_i)-s_{i+1}\ou^{i+1}+s_i\ou^i
\end{eqnarray}
and can further proceed in the following way that is illustrated by the examples below:\vspace*{-0.05in}

$\bullet$ If $\eta^i(t_i)>0$, it follows from (3) that $\la x^i_*,q(t_i)\ra=c_i$. Combining this with the maximization condition (6) allows us to determine an optimal control and the corresponding optimal motion dynamics.\vspace*{-0.05in}

$\bullet$ If $\eta^i(t_i)=0$, then the problem can be solved via \eqref{e:147}.\\[1ex]
Observe also that in our setting it is possible to represent the cost functional \eqref{cr-cost} as a function of $(\ou^1,\ldots,\ou^n)$ and $\eta^i(t_j)$ with $i=0,\ldots,n$ and $t_j\in[0,T]$. Thus the original optimal control problem can be reduced to finite-dimensional optimization of this cost subject to the constraints in \eqref{e:145} and \eqref{e:147}.\vspace*{-0.05in}

In the remainder of this section we consider two numerical examples with $n=2$ and $n=3$ participants, where the outlined procedure allows us to completely solve the formulated optimal control problem for the pedestrian traffic flow model. \vspace*{-0.1in}

\begin{example} {\bf(solving the controlled pedestrian traffic flow model with two participants).}\label{ex2} Specify the data of \eqref{e:135*} and \eqref{cr-cost} as follows:
$n=2,\;T=6,\;s_1=8,\;s_2=2,\;x^{01}=-60,\;x^{02}=-48,\;R=3$, and $c_j=-2R$ for $j=1,2$. Then the equations in \eqref{e:137*} reduce to
\begin{eqnarray}\label{e:139a}
\ox^1(t)=-60-\disp{\int^t_0}\eta(\tau)d\tau+ts_1\ou^1,
\end{eqnarray}\vspace*{-0.2in}
\begin{eqnarray}\label{e:139b}
\ox^2(t)=-48+\disp{\int^t_0}\eta(\tau)d\tau+ts_2\ou^2
\end{eqnarray}
for all $t\in[0,6]$. Define the convex and compact control set $U$ in  \eqref{e:135*} by
\begin{equation*}
U:=\big\{(u^1,u^2)\in\R^2\big|-1.8\le u^1=u^2\le 1.8\big\},
\end{equation*}
and let $t_1\in[0,6]$ be the first time when $\ox^2(t_1)-\ox^1(t_1)=2R=6$. If $t<t_1$, we get $x^2(t)-x^1(t)>2R=6$, and it follows from (2) that $\eta(t)=0$. At $t=t_1$ the motion $\ox(t)$ hits the state constraint set $C$ in \eqref{cr-state}, and hence it is reflected by a nonzero measure $\gg$ in (5). Now subtracting \eqref{e:139a} from \eqref{e:139b} with $t=t_1$ and taking into account that $\int_0^{t_1}\eta(\tau)d\tau=0$ tell us that
\begin{eqnarray}\label{oa}
12+t_1(2\ou^2-8\ou^1)=6,\;\mbox{ and so }\;-8\ou^1+2\ou^2+1\le 0\;\mbox{ by }\;t_1\le 6.
\end{eqnarray}
Suppose without loss of generality that both vector functions $\eta(t)$ and $\dot\ox(t)$ are well defined at $t=t_1$. Then we get from \eqref{e:135a} the equations
\begin{eqnarray}\label{velocities}
\left\{\begin{array}{ll}
\dot{\ox}^1(t_1)=-\eta(t_1)+8\ou^1,\\
\dot{\ox}^2(t_1)=\eta(t_1)+2\ou^2,
\end{array}\right.
\end{eqnarray}
which being combined with $\dot{\ox}^1(t_1)\le\dot{\ox}^2(t_1)$ give us the formulas
\begin{eqnarray}\label{etaex1}
-2\eta(t_1)+8\ou^1-2\ou^2\le 0.
\end{eqnarray}
Thus we deduce from \eqref{oa} and \eqref{etaex1} that $\eta(t_1)\ge 1/2$.\vspace*{-0.05in}

Remember that after the contact time $t_1$ both pedestrians tend to maintain their new constant velocities until $t=6$, and thus it holds that $\dot{\ox}(t)=\dot{\ox}(t_1)$ for all $t\in[t_1,6]$. Taking into account that $\ou(\cdot)$ is a constant on $[0,6]$ and that $\dot{\ox}(\cdot)$ is constant on the intervals $[0,t_1)$ and $[t_1,6]$, we get that the vector function $\eta(\cdot)$ is constant on $[0,t_1)$ and $[t_1,6]$, i.e.,
\begin{eqnarray}\label{eta_ex1}
\eta(t)=\left\{\begin{array}{ll}
\eta(0)&\mbox{ for a.e.}\;t\in[0,t_1)\;\mbox{ including }\;t=0,\\
\eta(t_1)&\mbox{ for a.e.}\;t\in[t_1,6]\;\mbox{ including }\;t=t_1.
\end{array}\right.
\end{eqnarray}
If $\eta(t)=\eta(t_1)>0$ a.e.\ on $[t_1,6]$, then it follows from (2) that $\ox^2(t)-\ox^1(t)=2R=6$ for all $t\in[t_1,6]$, and hence it shows that the optimal motion stays on the boundary of the state constraints \eqref{cr-state} on the entire interval $[t_1,6]$. Using further \eqref{e:139a}, \eqref{e:139b}, $\ox^2(t)-\ox^1(t)=6$ for all $t\in[t_1,6]$, \eqref{eta_ex1}, and the first equation in \eqref{oa}  gives us the relationships
\begin{eqnarray*}
&\;&x^2(t)-x^1(t)=12+2\(\int_0^{t_1}\eta(\tau)d\tau+\int_{t_1}^t\eta(\tau)d\tau\)+t\(s_2\ou^2-s_1\ou^1\)\\
&\Longleftrightarrow&6=12+2\(t-t_1\)\eta(t_1)+t\(s_2\ou^2-s_1\ou^1\)\\
&\Longleftrightarrow&12+t_1\(s_2\ou^2-s_1\ou^1\)=12+2\(t-t_1\)\eta(t_1)+t\(s_2\ou^2-s_1\ou^1\)\\
&\Longleftrightarrow&0=(t-t_1)\big[2\eta(t_1)-8\ou^1+2\ou^2\big],\quad t\in[t_1,6],
\end{eqnarray*}
which yield $2\eta(t_1)-s_1\ou^1+s_2\ou^2=0$. Combining this with the construction of the control set $U$ where $\ou^1=\ou^2$, we calculate the value of $\eta(\cdot)$ at the contact time $t=t_1$ by $\eta(t_1)=3\ou^2=3\ou^1$. Recalling that $\dot{\ox}^2(t_1)=\dot{\ox}^1(t_1)$ in this setting and remembering that $t_1=(\ou^2)^{-1}$ by \eqref{e:145} or \eqref{oa} and that $\eta(t_1)=3\ou^2$ allows us to express the value of cost functional \eqref{cr-cost} at $(\ox,\ou)$ by
\begin{eqnarray*}
J[\ox,\ou]&=&\frac{1}{2}\Big[\big(-60-(6-t_1)\eta(t_1)+6\cdot 8\ou^2\big)^2+\big(-48+(6-t_1)3\ou^2+6\cdot 2\ou^2\big)^2\Big]\\
&=&\frac{1}{2}\Big[\big(30\ou^2-57\big)^2+\big(30\ou^2-51\big)^2\Big].
\end{eqnarray*}
Minimizing the latter function of $\ou^2$ subject to the constraint $\ou^2\geq \frac{1}{6}$, which follows from the second inequality in \eqref{oa}, we get the optimal control value $\ou^2=\frac{3240}{1800}=1.8$.\vspace*{-0.05in}

Let us now calculate all the other elements of the optimal solution with the corresponding values of dual elements from the necessary optimality conditions. 
It follows from the first part of the maximization condition (6) that we can choose $\psi(t)=\ou(t)=\begin{pmatrix}
1.8\\
1.8
\end{pmatrix}$ and hence arrive at $q(t)=\begin{pmatrix}
q^{1}(t)\\
q^{2}(t)
\end{pmatrix}\equiv\begin{pmatrix}
0.225\\
0.9
\end{pmatrix}$ by the second part of (6) and then calculate
\begin{eqnarray*}
p(t)&=&\begin{pmatrix}
p^{1}(t)\\
p^{2}(t)
\end{pmatrix}=\begin{pmatrix}
p^{1}(6)\\
p^{2}(6)
\end{pmatrix}\\
&=&-\lm\ox(6)-\eta(6)x^1_*\\
&=&-\begin{pmatrix}
-60-(6-t_1)\cdot 3\ou^1+6\cdot 8\cdot\ou^1\\
-48+(6-t_1)\cdot 3\ou^2+6\cdot 2\ou^2
\end{pmatrix}-\eta(t_1)x^1_*\\
&=&-\begin{pmatrix}
-60-(6-t_1)\cdot 3\ou^1+6\cdot 8\cdot\ou^1\\
-48+(6-t_1)\cdot 3\ou^2+6\cdot 2\ou^2
\end{pmatrix}-3\begin{pmatrix}
1.8\\
-1.8
\end{pmatrix}\\
&=&\begin{pmatrix}
-2.4\\
2.4
\end{pmatrix}\;\mbox{ for all }\;t\in[0,6]
\end{eqnarray*}
due to (4) and (7). Then it follows from $\gamma\([t,6]\)=p(t)-q(t)$ by (5) that
\begin{eqnarray*}
\gg\big([t,6]\big)=\begin{pmatrix}
-2.6\\
1.5
\end{pmatrix}\;\mbox{ for }\;0.5556\approx t_1\le t\le 6,
\end{eqnarray*}
which shows that the optimal sweeping motion hits the boundary of the state constraints at $t_1\approx 0.5556$ and then stays there until $T=6$.
\end{example}\vspace*{-0.15in}

The next example concerns the case of three participants in the pedestrian traffic flow model.\vspace*{-0.1in}

\begin{example} {\bf (solving the controlled pedestrian traffic flow model with three participants).}\label{ex3} Consider the optimal control problem in \eqref{e:135*} and \eqref{cr-cost} with the following data: $n=3,\;s_1=8,\;s_2=4,\;s_3=2,\;T=6,\;R=3,\;x^{01}=-60,\;x^{02}=-48,\;x^{03}=-42,\;c_j=-2R$ for $j=1,2,3$, and the compact convex control set $U$ given by
\begin{equation*}
U:=\big\{(u^1,u^2,u^3)\in\R^3\big|\;\max\{|u^1|,\,|u^2|,\,|u^3|\}\le 2\big\}.
\end{equation*}
Following the procedure outlined above, we first obtain $x^{02}-x^{01}=12>6=2R$ and $x^{03}-x^{02}=6=2R$. Then it is obvious that $t_2=0$, and $t_1$ is calculated from \eqref{e:145} by
\begin{eqnarray*}
t_1=\dfrac{6}{\eta^2(0)-4\ou^2+8\ou^1}\le 6.
\end{eqnarray*}
Hence $\Theta_1=t_2=0$, and we get for $t\in[0,t_1)$ that
\begin{eqnarray*}
\dot{\ox}^1(t)=8\ou^1,\;\dot{\ox}^2(t)=-\eta^2(0)+4\ou^2,\;\dot{\ox}^3(t)=\eta^2(0)+2\ou^3,
\end{eqnarray*}\vspace*{-0.35in}
\begin{eqnarray*}
\left\{\begin{array}{ll}
\ox^1(t)=-60+8t\ou^1,\\
\ox^2(t)=-48-t\eta^2(0)+4t\ou^2,\\
\ox^3(t)=-42+t\eta^2(0)+2t\ou^3.
\end{array}\right.
\end{eqnarray*}
When $t\in[t_1,6]$, the corresponding representations of the velocities and the trajectories are
\begin{eqnarray*}
\dot{\ox}^1(t)=-\eta^1(t_1)+8\ou^1,\;\dot{\ox}^2(t)=-\eta^2(t_1)+\eta^1(t_1)+4\ou^2,\;\dot{\ox}^3(t)=\eta^2(t_1)+2\ou^3,
\end{eqnarray*}\vspace*{-0.25in}
\begin{eqnarray*}
\left\{\begin{array}{ll}
\ox^1(t)=-60-(t-t_1)\eta^1(t_1)+8t\ou^1,\\
\ox^2(t)=-48-t\eta^2(0)+(t-t_1)\big(\eta^1(t_1)-\eta^2(t_1)\big)+4t\ou^2,\\
\ox^3(t)=-42+t\eta^2(0)+(t-t_1)\eta^2(t_1)+2t\ou^3.
\end{array}\right.
\end{eqnarray*}
It follows from \eqref{e:147} and the obvious condition $\eta^1(0)=0$ that
\begin{eqnarray*}
2\eta^1(t_1)=\eta^2(t_1)-4\ou^2+8\ou^1,\;2\eta^2(0)=-2\ou^3+4\ou^2,\;\mbox{ and }\;2\eta^2(t_1)=\eta^1(t_1)-2\ou^3+4\ou^2.
\end{eqnarray*}
Using (3) together with $\eta^1(t_1)>0$ yields $\la x^1_*,q(t_1)\ra=c_1$ and hence $q^1(t_1)-q^2(t_1)=c_1=-6$. Then we can rewrite the above expressions for $\eta(\cdot)$ and the formula for $t_1$ as
\begin{eqnarray}\label{e:151}
\left\{\begin{array}{ll}
t_1=\disp\frac{6}{8\ou^1-2\ou^2-\ou^3},\\
\eta^2(0)=2\ou^2-\ou^3,\\
\eta^1(t_1)=\disp\frac{16}{3}\ou^1-\frac{4}{3}\ou^2-\frac{2}{3}\ou^3,\\
\eta^2(t_1)=\disp\frac{8}{3}\ou^1+\frac{4}{3}\ou^2-\frac{4}{3}\ou^3.
\end{array}\right.
\end{eqnarray}
Since pedestrians~2 and ~3 are in contact at the beginning and since we do not know whether $\eta^2(0)>0$ or $\eta^2(0)=0$, let consider the following two cases:\\[1ex]
{\bf Case~1: $\eta^2(0)>0$}. Then it follows from (3) that
\begin{eqnarray*}
\la x^2_*,q(0)\ra=c_2,\;\mbox{ i.e., }\;q^2(0)-q^3(0)=-6.
\end{eqnarray*}
Combining this with the equality $q^1(t_1)-q^2(t_1)=c_1=-6$ obtained above, we can choose
$q(t)=\begin{pmatrix}
1\\
7\\
13
\end{pmatrix}$
and then deduce from the formula for $\psi(t)$ in (6) that
\begin{eqnarray*}
\psi(t)=
\begin{pmatrix}
\psi^1(t)\\
\psi^2(t)\\
\psi^3(t)
\end{pmatrix}=\begin{pmatrix}
8q^1(t)\\
4q^2(t)\\
2q^3(t)
\end{pmatrix}\;\mbox{ on }\;[0,6].
\end{eqnarray*}
Thus $\psi(t)\equiv\begin{pmatrix}
8\\
28\\
26
\end{pmatrix}$ on $[0,6]$. Now the maximization condition (6) gives us the optimal control $\ou(t)=\begin{pmatrix}
2\\
2\\
2
\end{pmatrix}$, which lies on the boundary of the control set $U$. The corresponding optimal contact time and the optimal motion dynamics are, respectively, $t_1=0.6$ and
\begin{eqnarray*}
\big(\ox_1(t),\ox_2(t),\ox_3(t)\big)=\left\{\begin{array}{ll}
(16t-60,\;6t-48,\;6t-42)&\mbox{ for }\;t\in[0,t_1),\\
\(\disp\frac{28}{3}t-56,\;\frac{22}{3}t-48.8,\;\frac{34}{3}t-45.2\)&\mbox{ for }\;t\in[t_1,6].
\end{array}\right.
\end{eqnarray*}
Note also that $\gamma([t,6])=\begin{pmatrix}
-\frac{23}{3}\\
-\frac{13}{15}\\
-\frac{457}{15}
\end{pmatrix}$
when $t\in[t_1,6]$ for $\lm=1$. Using (5) and calculating $q(t)$ as above, we arrive then at the following calculation of $p(\cdot)$:
\begin{eqnarray*}
p(t)\equiv p(T)&=&-\ox(T)-\sum_{i\in I(\ox(T))}\eta^i(T)x^i_*\\
&=&-\begin{pmatrix}
\frac{28}{3}\cdot 6-56\\
\frac{22}{3}\cdot 6-48.8\\
\frac{34}{3}\cdot 6-45.2
\end{pmatrix}-\eta^1(t_1)\begin{pmatrix}
1\\
-1\\
0
\end{pmatrix}-\eta^2(t_1)\begin{pmatrix}
0\\
1\\
-1
\end{pmatrix}\\
&=&\begin{pmatrix}
-\frac{20}{3}\\
\frac{92}{15}\\
-\frac{262}{15}
\end{pmatrix}.
\end{eqnarray*}
{\bf Case~2:} $\eta^2(0)=0$. In this case we have $2\ou^2=\ou^3$, and it follows from \eqref{e:151} that
\begin{eqnarray*}
\left\{\begin{array}{ll}
t_1=\disp\frac{6}{8\ou^1-4\ou^2}\le 6,\\
\eta^2(0)=0,\\
\eta^1(t_1)=\frac{16}{3}\ou^1-\frac{8}{3}\ou^2,\\
\eta^2(t_1)=\frac{8}{3}\ou^1-\frac{4}{3}\ou^2.
\end{array}\right.
\end{eqnarray*}
Since $\eta^1(t_1)>0$ and $\eta^2(t_1)>0$, we get $2\ou^1>\ou^2$, and hence we can choose $q(t)=\begin{pmatrix}
1\\
7\\
13
\end{pmatrix}$
by (3). The maximization condition (6) gives us $\ou^1=\ou^2=\ou^3=2$, which contradicts the relationships $2\ou^1>\ou^2=\frac{\ou^3}{2}$ obtained above. This tells us that the situation in Case~2 cannot be realized, and therefore the calculations in Case~1 completely solve the problem under consideration in this example.
\end{example}
\vspace*{-0.4in}

\section{Conclusions}\vspace*{-0.15in}

This paper presents applications of recently obtained results on optimal control of perturbed sweeping processes to two practical models known as the mobile robot model with obstacle and the pedestrian traffic flow model through a doorway. We see that the approach and developments of \cite{cmn}, based on advanced variational analysis and the method discrete approximations, provide efficient tools to determine optimal solutions to naturally formulated control versions of these models via new necessary optimality conditions expressed entirely in terms of the model data. Nontrivial numerical examples presented in the paper give us exact solutions of the control problems formulated for the models under consideration in the case of lower numbers of participants and illustrate the scheme of applications of the obtained necessary optimality conditions in more general settings.\vspace*{-0.05in}

Our further research goals concerning these models include developing efficient numerical algorithms to solve the optimal control problems for them with large numbers of participants. It could be done, in particular, by using an appropriate discretization and employing numerical algorithms of finite-dimensional optimization to the discrete-time problems obtained in this way. We also believe that the developed necessary optimality conditions for the perturbed sweeping processes would be useful to investigate other practical model with a sweeping process dynamics that frequently appear in various branches of mechanics, engineering, economics, behavioral sciences, etc. \vspace*{-0.1in}

\begin{acknowledgements}
Research of the second and third authors was partly supported by the USA National Science Foundation under grants DMS-1512846 and DMS-1808978, and by the USA Air Force Office of Scientific Research under grant \#15RT0462. The authors are thankful to Tan Cao for many useful discussions.
\end{acknowledgements}
\end{document}